# LARGE DEVIATIONS OF THE EMPIRICAL VOLUME FRACTION FOR STATIONARY POISSON GRAIN MODELS

By Lothar Heinrich

*University of Augsburg*

We study the existence of the (thermodynamic) limit of the scaled cumulant-generating function $L_n(z) = |W_n|^{-1} \log \mathsf{E} \exp\{z|\Xi \cap W_n|\}$ of the empirical volume fraction $|\Xi \cap W_n|/|W_n|$, where $|\cdot|$ denotes the $d$-dimensional Lebesgue measure. Here $\Xi = \bigcup_{i \geq 1}(\Xi_i + X_i)$ denotes a $d$-dimensional Poisson grain model (also known as a Boolean model) defined by a stationary Poisson process $\Pi_\lambda = \sum_{i \geq 1} \delta_{X_i}$ with intensity $\lambda > 0$ and a sequence of independent copies $\Xi_1, \Xi_2, \ldots$ of a random compact set $\Xi_0$. For an increasing family of compact convex sets $\{W_n, n \geq 1\}$ which expand unboundedly in all directions, we prove the existence and analyticity of the limit $\lim_{n \to \infty} L_n(z)$ on some disk in the complex plane whenever $\mathsf{E} \exp\{a|\Xi_0|\} < \infty$ for some $a > 0$. Moreover, closely connected with this result, we obtain exponential inequalities and the exact asymptotics for the large deviation probabilities of the empirical volume fraction in the sense of Cramér and Chernoff.

**1. Introduction and main results.** The *Poisson grain model* (PGM; also known as the *Boolean model*) is the best studied and most used random set model to describe systems of randomly distributed and irregularly shaped clumps in a Euclidean space $R^d$, $d \geq 1$ [see Matheron (1975), Hall (1988) or Stoyan, Kendall and Mecke (1995)]. It is the basic model in stereology and stochastic geometry. Statistical analysis of a stationary PGM is mostly based on a single realization of the union set of clumps in some region $W$ which is assumed to expand unboundedly in all directions [see, e.g., Molchanov (1997)]. To be definite in describing our problem, we first give a rigorous definition of a stationary PGM as the union set

$$\Xi := \bigcup_{i \geq 1} (\Xi_i + X_i) \tag{1.1}$$

of independent copies $\Xi_1, \Xi_2, \ldots$ (grains) of a random compact set $\Xi_0$ (typical grain) that has distribution $Q$, where the grains are independently shifted by the atoms $X_1, X_2, \ldots$ (germ points) of a stationary Poisson process $\Pi_\lambda = \sum_{i \geq 1} \delta_{X_i}$ with intensity $\lambda$ (= mean number of germ points in the unit cube $[0,1)^d$). Throughout this paper, all random elements are defined on a common probability space $[\Omega, \mathfrak{A}, \mathsf{P}]$ and $\mathsf{E}$ denotes the expectation with respect to $\mathsf{P}$. In particular, $\Xi_0$ is a measurable mapping from $[\Omega, \mathfrak{A}, \mathsf{P}]$ into the space of nonvoid compact subsets $\mathcal{K}$ of $R^d$ equipped with the Hausdorff metric and $Q$ coincides with the image measure $\mathsf{P} \circ \Xi_0^{-1}$ that acts on the corresponding Borel $\sigma$-field $\mathfrak{B}(\mathcal{K})$ [see Matheron (1975)]. Note that $\Xi$ is a closed set (P-a.s.) if $\mathsf{E}|\Xi_0 + B_r(o)| < \infty$ for $r > 0$, where $B_r(x)$ denotes the closed ball with radius $r > 0$ centered at $x \in R^d$ and $|\cdot|$ denotes the Lebesgue measure in $R^d$ [see Heinrich (1992)].

The main aim of this paper is to prove the existence and analyticity of the limit (as $n \to \infty$) of

$$L_n(z) := \frac{1}{|W_n|} \log \mathsf{E} \exp\{z|\Xi \cap W_n|\} \qquad \text{on } D_\Delta := \{z \in C^1 : |z| < 1/\Delta\} \tag{1.2}$$

for some $0 < \Delta < \infty$ provided that an exponential moment of the volume $|\Xi_0|$ exists, that is,

$$M(a) := \mathsf{E} \exp\{a|\Xi_0|\} < \infty \qquad \text{for some } a > 0, \tag{1.3}$$

and $(W_n)$ is a *convex averaging sequence* of sets in $R^d$, that is, each $W_n$ is a (deterministic) compact convex set, $(W_n)$ is nondecreasing and its union is $R^d$ [see Daley and Vere-Jones (1988)]. Because of the conspicuous analogy to similar problems in statistical physics [see Ruelle (1969)], we call $L(z) = \lim_{n \to \infty} L_n(z)$ the *thermodynamic limit* of (the thermodynamic function) $L_n(z)$. The second aim, which is closely connected with the first, is to derive inequalities and asymptotic relationships (in the sense of Cramér and Chernoff) for probabilities of large deviations of the empirical volume fraction $\hat{p}_n := |\Xi \cap W_n|/|W_n|$ from its mean $p := \mathsf{E}|\Xi \cap [0,1)^d| = \mathsf{P}(o \in \Xi)$.

In the special case of a bounded typical grain, that is, $\Xi_0 \subseteq B_R(o)$ for some $0 < R < \infty$, both problems were solved satisfactorily by Götze, Heinrich and Hipp (1995) using the device of $m$-dependent random fields with block representation. The proving technique in the present paper is completely different from that in Götze, Heinrich and Hipp (1995) and does not require any strong mixing properties of the PGM (1.1) as one would expect. In general, (1.3) does not imply specific mixing rates as needed, for example, in Mase (1982) or Heinrich and Molchanov (1999). However, in case of a spherical typical grain, (1.3) induces an exponentially decaying $\beta$-mixing coefficient [see Heinrich and Molchanov (1999)]. Note that (1.3) does not even imply the closedness of $\Xi$ in general; see the Appendix. For a positive random variable $X$ with infinite mean, the typical grain $\Xi_0 = [0, X] \times [0, 1/X]$ exhibits such an example for $d = 2$.

For this reason we choose the probability space $[\Omega, \mathfrak{A}, \mathsf{P}]$ (for its existence, see the Appendix) in such a way that the mapping $R^d \times \Omega \ni (x, \omega) \mapsto \mathbf{1}_{\Xi(\omega)}(x)$ is $\mathfrak{B}(R^d) \otimes \mathfrak{A}$-measurable. This property allows us to apply Fubini's theorem to the 0–1-valued random field $\xi(x) = \mathbf{1}_\Xi(x)$, $x \in R^d$, and implies that the function

$$(1.4) \qquad p_\Xi^{(k)}(x_1, \ldots, x_k) := \mathsf{E}\xi(x_1) \cdots \xi(x_k) = \mathsf{P}(x_1 \in \Xi, \ldots, x_k \in \Xi)$$

is $\mathfrak{B}(R^{dk})$-measurable for any $k \geq 1$ and $|\Xi \cap W| = \int_W \xi(x)\, dx$ is a random variable over $[\Omega, \mathfrak{A}, \mathsf{P}]$ for any bounded $W \in \mathfrak{B}(R^d)$. The functions (1.4) are expressible (and vice versa) by the corresponding probabilities for the complement set $\Xi^c$:

$$\begin{aligned}(1.5) \qquad p_{\Xi^c}^{(k)}(x_1, \ldots, x_k) &:= \mathsf{E}(1 - \xi(x_1)) \cdots (1 - \xi(x_k)) \\ &= \mathsf{P}(\Xi \cap \{x_1, \ldots, x_k\} = \varnothing).\end{aligned}$$

Since $(\Xi_i + X_i) \cap \{x_1, \ldots, x_k\} = \varnothing$ iff $X_i \notin (-\Xi_i) + \{x_1, \ldots, x_k\}$, the shape of the probability generating functional (A.2) (of a stationary independently marked Poisson process $\Pi_{\lambda, Q}$) for $v(x, K) = 1 - \mathbf{1}_{(-K)+\{x_1, \ldots, x_k\}}(x)$ yields

$$\begin{aligned}(1.6) \qquad p_{\Xi^c}^{(k)}(x_1, \ldots, x_k) &= \mathsf{E}\prod_{i \geq 1}(1 - \mathbf{1}_{(-\Xi_i)+\{x_1, \ldots, x_k\}}(X_i)) \\ &= \exp\left\{-\lambda\left|\bigcup_{i=1}^k (\Xi_0 - x_i)\right|\right\}.\end{aligned}$$

Note that $p_{\Xi^c}^{(k)}(x_1, \ldots, x_k) = 1 - T_\Xi(\{x_1, \ldots, x_k\})$ for an arbitrary random closed set $\Xi$ with capacity functional $T_\Xi$ [see Matheron (1975)]. The study of the sequence (1.2) is closely related with the behavior of the higher-order *mixed cumulants*

$$(1.7) \qquad c_\Xi^{(k)}(x_1, \ldots, x_k) := \Gamma(\xi(x_1), \ldots, \xi(x_k)) \qquad \text{for } k \geq 1$$

of the random field $\{\xi(x), x \in R^d\}$, where the *mixed cumulant (semi-invariant)* of any random variables $Y_1, \ldots, Y_k$ (having a finite $k$th moment) is defined by

$$(1.8) \quad \Gamma(Y_1, \ldots, Y_k) := i^{-k} \frac{\partial^k}{\partial s_1 \cdots \partial s_k} \log \mathsf{E}\exp\left\{i\sum_{j=1}^k s_j Y_j\right\}\bigg|_{s_1 = \cdots = s_k = 0}$$

and $\Gamma_k(Y) := \Gamma(Y,\ldots,Y)$ [obtained by putting $Y = Y_1 = \cdots = Y_k$ in (1.8)] denotes the $k$th *cumulant* of $Y$. Directly from (1.8) it is seen that, for $k \geq 2$,

$$(1.9) \quad c^{(k)}_{\Xi^c}(x_1,\ldots,x_k) := \Gamma(1-\xi(x_1),\ldots,1-\xi(x_k)) = (-1)^k c^{(k)}_{\Xi}(x_1,\ldots,x_k).$$

We are now in a position to formulate our main result.

THEOREM 1. *Let $\Xi$ be the PGM (1.1) with compact typical grain $\Xi_0$ satisfying (1.3) and let $\{W_n, n \geq 1\}$ be a convex averaging sequence in $R^d$. Then, for any $k \geq 2$,*

$$(1.10) \quad \int_{(R^d)^{k-1}} |c^{(k)}_{\Xi}(o,x_2,\ldots,x_k)|\,d(x_2,\ldots,x_k) \leq (k-1)! H(a) \Delta(a)^{k-2},$$

*where $H(a) := 8\lambda M(a)(1+\exp\{\lambda\mathsf{E}|\Xi_0|\})/a^2$ and $\Delta(a) := 8(a+\lambda M(a)) \times (1+\exp\{\lambda\mathsf{E}|\Xi_0|\})/a^2$. Furthermore, the limit $L(z) = \lim_{n\to\infty} L_n(z)$ exists and is analytic on the open disk $D_{\Delta(a)}$.*

The next result states Cramér's large deviations relationships for the random sequence $|\Xi \cap W_n|$ and an optimal Berry–Esseen bound of the distance between $F_n(x) := \mathsf{P}(\sqrt{|W_n|}(\hat{p}_n - p) \leq x\sigma_n)$ and the standard normal distribution function $\Phi(x) = \int_{-\infty}^x \exp(-t^2/2)\,dt/\sqrt{2\pi}$, where

$$\sigma_n^2 := \frac{\mathsf{Var}(|\Xi \cap W_n|)}{|W_n|}$$
$$= \int_{R^d} \frac{|W_n \cap (W_n - x)|}{|W_n|}(\exp(-\lambda|\Xi_0 \cup (\Xi_0 - x)|) - \exp(-2\lambda|\Xi_0|))\,dx.$$

The following Theorem 2 is derived from (1.10) combined with a well-known lemma on large deviations for a single random variable discussed by Statulevičius (1966) [see also Saulis and Statulevičius (1991), Lemma 2.3].

THEOREM 2. *Let the assumptions of Theorem 1 be satisfied and, in addition, let $\mathsf{E}|\Xi_0| > 0$. Then $\sigma_n^2$ converges to a nonzero limit $\int_{R^d} c^{(2)}_{\Xi}(o,x)\,dx$ and, for $0 \leq x \leq \sigma_n\sqrt{|W_n|}/2\Delta(a)(1+4H_n)$ with $H_n = H(a)/2\sigma_n^2$, the asymptotic relationships*

$$(1.11) \quad \frac{1-F_n(x)}{1-\Phi(x)} = \exp\left\{\frac{x^3}{\sigma_n\sqrt{|W_n|}} \sum_{k=0}^\infty \mu_k^{(n)}\left(\frac{x}{\sigma_n\sqrt{|W_n|}}\right)^k\right\}\left(1+O\left(\frac{1+x}{\sqrt{|W_n|}}\right)\right)$$

*and*

$$(1.12) \quad \frac{F_n(-x)}{\Phi(-x)} = \exp\left\{\frac{-x^3}{\sigma_n\sqrt{|W_n|}} \sum_{k=0}^\infty \mu_k^{(n)}\left(\frac{-x}{\sigma_n\sqrt{|W_n|}}\right)^k\right\}\left(1+O\left(\frac{1+x}{\sqrt{|W_n|}}\right)\right)$$

*hold as* $n \longrightarrow \infty$, *where the coefficients*

$$\mu_k^{(n)} = \frac{1}{(k+2)(k+3)}$$

(1.13)
$$\times \sum_{l=1}^{k+1}(-1)^{l-1}\binom{k+l+1}{l}\sum_{\substack{k_1+\cdots+k_l=k+1 \\ k_i\geq 1,\ i=1,\ldots,l}}\prod_{i=1}^{l}\frac{\Gamma_{k_i+2}(|\Xi\cap W_n|)}{\sigma_n^2|W_n|(k_i+1)!}$$

*satisfy the estimate* $|\mu_k^{(n)}| \leq 4H_n\Delta(a)(2\Delta(a)(1+4H_n))^k/(k+2)(k+3)$ *for* $k \geq 0$. *Furthermore, there exists some constant* $c > 0$ [*depending on* $a$, $\lambda$, $M(a)$ *and* $\sigma_n^2$] *such that*

$$(1.14) \qquad \sup_{x\in R^1}|F_n(x) - \Phi(x)| \leq \frac{c}{\sqrt{|W_n|}}.$$

Our next Theorem 3 provides large deviations inequalities for the unbiased estimators

$$(1.15)\quad \hat{p}_W := \frac{|\Xi\cap W|}{|W|} \quad \text{and}\quad \hat{C}_W(x) := \frac{|\Xi\cap(\Xi-x)\cap W|}{|W|},\qquad x\in R^d,$$

of the volume fraction $p = \mathsf{P}(o\in\Xi)$ and the covariance $C(x) = \mathsf{P}(o\in\Xi, x\in\Xi)$, respectively, in the case when the PGM (1.1) is observed on a sampling window $W \in \mathfrak{B}(R^d)$. Note that, in contrast to the volume fraction $p$, the covariance $C(\cdot)$ reveals information on the inner structure of the random set $\Xi$ [see Matheron (1975) and Stoyan, Kendall and Mecke (1995)]. The deviation of the estimators (1.15) from their means $p$ and $C(x)$ is estimated under finite-order as well as exponential moment assumptions put on the volume of the typical grain $\Xi_0$.

THEOREM 3. *Let* $\Xi$ *be the PGM* (1.1) *with compact typical grain* $\Xi_0$ *that satisfy* $\mathsf{E}|\Xi_0|^s < \infty$ *for some real* $s \geq 2$. *Furthermore, let* $W \subset R^d$ *be a bounded Borel set with inner points. Then there exist positive constants* $c_s^{(1)}(\lambda)$ *and* $c_s^{(2)}(\lambda)$ (*depending on* $\lambda$ *and the moments* $\mathsf{E}|\Xi_0|^k$, $k = 1,\ldots,[s], s$) *such that, for any* $\varepsilon > 0$,

$$(1.16) \qquad \mathsf{P}(|\hat{p}_W - p| \geq \varepsilon) \leq c_s^{(1)}(\lambda)\varepsilon^{-s}|W|^{-s/2}$$

*and*

$$(1.17)\ \mathsf{P}(|\hat{C}_W(x) - C(x)| \geq \varepsilon) \leq c_s^{(2)}(\lambda)\varepsilon^{-s}|W|^{-s/2} \qquad \text{for all } x\in R^d.$$

*If* $\Xi_0$ *satisfies condition* (1.3), *then the Bernstein-type inequality*

$$(1.18)\quad \mathsf{P}(\hat{p}_W - p \geq \varepsilon) \leq \begin{cases} \exp\left\{-\dfrac{1-\rho}{2H(a)}\varepsilon^2|W|\right\}, & \text{if } 0 \leq \varepsilon \leq \dfrac{H(a)\rho}{\Delta(a)(1-\rho)}, \\ \exp\left\{-\dfrac{\rho}{2\Delta(a)}\varepsilon|W|\right\}, & \text{if } \varepsilon \geq \dfrac{H(a)\rho}{\Delta(a)(1-\rho)}, \end{cases}$$

*holds for any $0 < \rho < 1$ and $H(a), \Delta(a)$ from Theorem* 1. *Exactly the same bounds hold for the probability* $\mathsf{P}(\hat{p}_W - p \leq -\varepsilon)$.

In Theorem 4 below, we derive a Chernoff rate function [see Dembo and Zeitouni (1998) and references therein] for the sequence of empirical volume fractions $\hat{p}_n$ in terms of the thermodynamic limit $L(z)$, which provides an extension and refinement of the relationship (1.12) for the $x$ values $x(\varepsilon) = \varepsilon\sqrt{|W_n|}/\sigma_n$ with $\varepsilon \in (0, \varepsilon^*)$, where $\varepsilon^*$ is determined by the slope of the function $L(h)$ at $h = 1/\Delta(a)$.

THEOREM 4. *Under the assumptions of Theorem* 2 *the large deviations relationship*

$$(1.19) \qquad \lim_{n\to\infty} \frac{\log \mathsf{P}(\hat{p}_n - p \geq \varepsilon)}{|W_n|} = \inf_{0 \leq h < 1/\Delta(a)} g(h) = g(h_0(\varepsilon))$$

*holds in the interval* $0 \leq \varepsilon < \varepsilon^* := \lim_{h\uparrow 1/\Delta(a)} L'(h) - p$, *where* $g(h) := L(h) - h(\varepsilon + p)$ *and* $h_0(\varepsilon)$ *is the unique solution of the equation* $g'(h) = 0$, *that is*, $L'(h) = \varepsilon + p$. *A corresponding relationship is valid for the probability* $\mathsf{P}(\hat{p}_n - p \leq -\varepsilon)$, *where the function $g(h)$ is defined for $h \in (-1/\Delta(a), 0]$ and with $-\varepsilon$ instead of $\varepsilon$.*

This result touches the question of whether $\hat{p}_n$ satisfies the large deviation principle, the answer to which seems to be unknown so far. Without giving details we mention only that the limit $\lim_{n\to\infty} L_n(h)$ exists, on the negative real axis, which can be shown by the methods of Ruelle [(1969), Chapter 3.4]. For related problems concerning large deviation principles for stationary independently marked Poisson processes, refer to Georgii and Zessin (1993). Similar results for Young measures related to Poisson grain models have been proved by Piau (1999).

The rest of this paper is organized as follows: In Section 2, we investigate (1.2) for a quite general random set model (Lemmas 1 and 2) and put together the required tools from point process theory presented in a rather general setting (Lemma 3). In Section 3, we are concerned with the proof of Theorem 1, which is divided into several steps (Lemmas 4–7), whereas the proofs of the Theorems 2, 3 and 4 are deferred to Section 4. The Appendix contains, among other things, the construction of a measurable random field $\xi(x) = \mathbf{1}_\Xi(x)$, $x \in R^d$, and a criterion for (non-)closedness of the PGM $\Xi$ given by (1.1).

**2. Preliminary results and relationships to point processes.** We first investigate the behavior of the cumulants $\Gamma_k(|\Xi \cap W_n|)$ and give a condition which guarantees the existence and analyticity of the limit of (1.2) for the support set $\Xi = \mathrm{supp}(\xi)$ of an arbitrary $(\mathfrak{B}(R^d) \otimes \mathfrak{A})$-measurable, 0–1-valued, stationary random field $\{\xi(x), x \in R^d\}$. We use the same notation as in Section 1. Lemma 2 states that this condition can be expressed by the total variation of the reduced cumulant measures of the Cox process

$$(2.1) \qquad \Pi_{\Xi}^{(z)} := \sum_{i \geq 1} (1 - \mathbf{1}_{\Xi}(Y_i)) \delta_{Y_i},$$

which is directed by the random measure $z \int_{R^d} \mathbf{1}_{(\cdot)}(x)(1 - \xi(x))\,dx = z|\Xi^c \cap (\cdot)|$, where $\Pi_z = \sum_{i \geq 1} \delta_{Y_i}$ is a stationary Poisson process with intensity $z > 0$ that is independent of $\Xi$. In the second part of this section we introduce a family of correlation measures for arbitrary stationary point processes and derive (Lemma 3) a recurrence relationship for the corresponding Lebesgue density functions provided they exist. Lemma 3 is the key to prove Theorem 1 and it seems to be of interest on its own.

LEMMA 1. *Let $\{\xi(x), x \in R^d\}$ be a measurable, 0–1-valued, stationary random field on $R^d$ with support set $\Xi := \{x \in R^d : \xi(x) = 1\}$. Then, for any bounded $W \in \mathfrak{B}(R^d)$ and $k \geq 2$, we have*

$$|\Gamma_k(|\Xi \cap W|)| \leq |W| G_k(\Xi)$$

$$\text{with } G_k(\Xi) := \int_{(R^d)^{k-1}} |c_{\Xi}^{(k)}(o, x_2, \ldots, x_k)|\, d(x_2, \ldots, x_k).$$

*Furthermore, let $(W_n)_{n \geq 1}$ be a convex averaging sequence. If $G_k(\Xi) < \infty$ for some $k \geq 2$, then*

$$(2.2) \quad \lim_{n \to \infty} \frac{\Gamma_k(|\Xi \cap W_n|)}{|W_n|} = \int_{(R^d)^{k-1}} c_{\Xi}^{(k)}(o, x_2, \ldots, x_k)\, d(x_2, \ldots, x_k),$$

*and if $G_k(\Xi) \leq k! H \Delta^{k-2}$ for some $H, \Delta > 0$ and any $k \geq 2$, then the thermodynamic limit $L(z) = \lim_{n \to \infty} L_n(z)$ exists and is an analytic function on the open disk $D_\Delta$, where $L_n(z)$ is given by (1.2) with $\Xi := \{x \in R^d : \xi(x) = 1\}$ [instead of (1.1)]. For $z \in D_\Delta$, the function $L(z)$ admits the power series expansion*

$$L(z) = zp + \frac{z^2}{2} \int_{R^d} (C(x) - p^2)\, dx$$
$$+ \sum_{k \geq 3} \frac{z^k}{k!} \int_{(R^d)^{k-1}} c_{\Xi}^{(k)}(o, x_2, \ldots, x_k)\, d(x_2, \ldots, x_k),$$

*where $p := p_{\Xi}^{(1)}(o)$ (volume fraction of $\Xi$) and $C(x) := p_{\Xi}^{(2)}(o, x)$ (covariance of $\Xi$).*

PROOF. Using Fubini's theorem and the definition (1.4) we may write

$$
\begin{aligned}
(2.3) \quad \mathsf{E}\prod_{i=1}^{k}|\Xi\cap B_i|^k &= \mathsf{E}\prod_{i=1}^{k}\int_{B_i}\xi(x_i)\,dx_i \\
&= \int_{B_1\times\cdots\times B_k} p_{\Xi}^{(k)}(x_1,\ldots,x_k)\,d(x_1,\ldots,x_k).
\end{aligned}
$$

A direct calculation of the logarithmic derivatives in (1.8) leads to

$$
(2.4)\quad \Gamma(Y_1,\ldots,Y_k)=\sum_{j=1}^{k}(-1)^{j-1}(j-1)!\sum_{K_1\cup\cdots\cup K_j=K}\prod_{i=1}^{j}\mathsf{E}\left(\prod_{k_i\in K_i}Y_{k_i}\right)
$$

[see, e.g., Saulis and Statulevičius (1991)], where the inner sum is taken over all decompositions of $K=\{1,\ldots,k\}$ into $j$ disjoint nonempty subsets $K_1,\ldots,K_j$. From (2.4) and (2.3) and by repeated application of Fubini's theorem, we see that the integral $\int_{B_1\times\cdots\times B_k}\Gamma(\xi(x_1),\ldots,\xi(x_k))\,d(x_1,\ldots,x_k)$ coincides with $\Gamma(|\Xi\cap B_1|,\ldots,|\Xi\cap B_k|)$. This means, setting $B_1=\cdots=B_k=W$ and using (1.7), that

$$
(2.5)\qquad \Gamma_k(|\Xi\cap W|)=\int_{W^k}c_{\Xi}^{(k)}(x_1,x_2,\ldots,x_k)\,d(x_1,x_2,\ldots,x_k).
$$

The stationarity of the random field $\{\xi(x),x\in R^d\}$ implies the invariance of the mixed cumulants (1.7) under diagonal shifts, that is,

$$
(2.6)\qquad c_{\Xi}^{(k)}(x_1,x_2,\ldots,x_k)=c_{\Xi}^{(k)}(o,x_2-x_1,\ldots,x_k-x_1),
$$

whence, by substituting $y_j=x_j-x_1$, $j=2,\ldots,k$, it follows that

$$
\begin{aligned}
\Gamma_k(|\Xi\cap W|) &= \int_{(R^d)^{k-1}}\int_{R^d}c_{\Xi}^{(k)}(o,y_2,\ldots,y_k)\mathbf{1}_W(x_1) \\
&\qquad\times\prod_{j=2}^{k}\mathbf{1}_W(y_j-x_1)\,dx_1\,d(y_2,\ldots,y_k) \\
&= \int_{(R^d)^{k-1}}c_{\Xi}^{(k)}(o,x_2,\ldots,x_k) \\
&\qquad\times|W\cap(W-x_2)\cap\cdots\cap(W-x_k)|\,d(x_2,\ldots,x_k),
\end{aligned}
$$

proving the first part of Lemma 1. The limit (2.2) is an immediate consequence of Lebesgue's dominated convergence theorem and the fact that, in view of the geometric properties of the $W_n$'s [see Fritz (1970)],

$$
\lim_{n\to\infty}\frac{|W_n\cap(W_n-x_1)\cap\cdots\cap(W_n-x_{k-1})|}{|W_n|}=1
$$

for any fixed $x_1,\ldots,x_{k-1}\in R^d$.

The power series expansion of (1.2) is

$$L_n(z) = pz + \sum_{k \geq 2} \frac{\Gamma_k(|\Xi \cap W_n|)}{|W_n|} \frac{z^k}{k!}$$

and, hence, by our assumptions,

$$|L_n(z) - pz| \leq \sum_{k \geq 2} G_k(\Xi) \frac{|z|^k}{k!}$$

$$\leq |z|^2 H \sum_{k \geq 2} (|z|\Delta)^{k-2} = \frac{|z|^2 H}{1 - |z|\Delta} \qquad \text{for } |z| < \frac{1}{\Delta}.$$

Thus, for any $n \geq 1$, $L_n(z)$ is analytic on the open disk $D_\Delta$ and, by (2.2), $L_n(z)$ converges to $L(z)$ uniformly in any closed subset of $D_\Delta$, proving the analyticity of $L(z)$ on $D_\Delta$. $\square$

To obtain estimates of the form $G_n(\Xi) \leq n! H \Delta^{n-2}$ in the case of the PGM (1.1), we first show that $(-z)^n c_\Xi^{(n)}$ coincides with the $n$th-order cumulant density of the Cox process (2.1). In the second step we introduce a family of correlation measures $\gamma_\Psi^{(m,n)}$ and their Lebesgue densities, and study them for $\Psi = \Pi_\Xi^{(1)}$. In Section 3 we perform a somewhat involved and rather lengthy inductive estimation technique to derive bounds of the total variation of these correlation measures in terms of moments of $|\Xi_0|$ when $\Xi$ is given by (1.1). The basic idea of this method goes back to Ruelle (1964) [see also Ruelle (1969), Chapter 4.4], who developed it (without using the terminology of point processes) to prove the existence of thermodynamic limits for grand canonical Gibbs ensembles with pair interactions. An extension to ensembles with higher-order interactions was tried by Greenberg (1971), but it fails in our situation.

To begin with, we briefly recall the definition of the $n$th-order factorial moment (and cumulant) measure $\alpha_\Psi^{(n)}$ (and $\gamma_\Psi^{(n)}$) of a point process $\Psi = \sum_{i \geq 1} \delta_{Z_i}$ that satisfies $\mathsf{E}\Psi^n(K) < \infty$ for $K \in \mathcal{K}$ by means of its probability generating functional $G_\Psi[w] = \mathsf{E}(\prod_{i \geq 1} w(Z_i))$, where $w : R^d \mapsto [0,1]$ is Borel measurable such that $1 - w$ has bounded support [see, e.g., Daley and Vere-Jones (1988)]. Setting $w_{v_1,\ldots,v_n}^{B_1,\ldots,B_n}(x) = 1 + \sum_{j=1}^n (v_j - 1)\mathbf{1}_{B_j}(x)$ for $1 - \frac{1}{n} \leq v_j \leq 1$, $i = 1, \ldots, n$, and bounded $B_1, \ldots, B_n \in \mathfrak{B}(R^d)$, we define

$$\alpha_\Psi^{(n)}\left(\underset{j=1}{\overset{n}{\times}} B_j\right) := \lim_{v_1,\ldots,v_n \uparrow 1} \frac{\partial^n}{\partial v_1 \cdots \partial v_n} G_\Psi[w_{v_1,\ldots,v_n}^{B_1,\ldots,B_n}]$$

and

$$\gamma_\Psi^{(n)}\left(\underset{j=1}{\overset{n}{\times}} B_j\right) := \lim_{v_1,\ldots,v_n \uparrow 1} \frac{\partial^n}{\partial v_1 \cdots \partial v_n} \log G_\Psi[w_{v_1,\ldots,v_n}^{B_1,\ldots,B_n}].$$

If $\alpha_\Psi^{(n)}$ (resp. $\gamma_\Psi^{(n)}$) is absolutely continuous with respect to the Lebesgue measure on $R^{dn}$, then we denote the corresponding (factorial) moment (resp. cumulant density) by $p_\Psi^{(n)}$ (resp. $c_\Psi^{(n)}$). In the sequel we often write $p_\Psi^{(n)}(X_n)$ instead $p_\Psi^{(n)}(x_1,\ldots,x_n)$, where $X_n$ stands for the (unordered) point set $\{x_1,\ldots,x_n\}$. In case the point process $\Psi$ is stationary, there exists a unique (signed) measure $\gamma_{\Psi,\mathrm{red}}^{(n)}$ on $\mathfrak{B}(R^{d(n-1)})$—called $n$th-order *reduced cumulant measure*—such that

$$\gamma_\Psi^{(n)}\left(\underset{j=1}{\overset{n}{\times}} B_j\right) = \int_{B_1} \gamma_{\Psi,\mathrm{red}}^{(n)}\left(\underset{j=2}{\overset{n}{\times}} (B_j - x)\right) dx$$

(2.7)
$$\text{for any bounded } B_1,\ldots,B_n \in \mathfrak{B}(R^d).$$

Finally, a stationary point process $\Psi$ is said to be *Brillinger-mixing* [see, e.g., Ivanoff (1982) or Heinrich and Schmidt (1985)], if $\mathsf{E}\Psi^n([0,1)^d) < \infty$ and the total variation $\mathrm{var}(\gamma_{\Psi,\mathrm{red}}^{(n)})$ on $R^{d(n-1)}$ is finite for all $n \geq 2$.

LEMMA 2. *Let $\Xi$ be the support set of a measurable, 0–1-valued, stationary random field $\{\xi(x), x \in R^d\}$. Then the $n$th-order reduced cumulant measure $\gamma_{\Pi_\Xi^{(z)},\mathrm{red}}^{(n)}$ of the Cox process (2.1) exists for any $n \geq 2$ and its total variation (if it exists) takes the form*

$$\mathrm{var}(\gamma_{\Pi_\Xi^{(z)},\mathrm{red}}^{(n)}) = z^n G_n(\Xi)$$

$$\text{with } \gamma_{\Pi_\Xi^{(z)},\mathrm{red}}^{(n)}(B) = (-z)^n \int_B c_\Xi^{(n)}(o, x_2, \ldots, x_n)\, d(x_2,\ldots,x_n)$$

*for $B \in \mathfrak{B}(R^{d(n-1)})$. Consequently, $\Pi_\Xi^{(z)}$ is Brillinger-mixing iff $G_n(\Xi) < \infty$ for all $n \geq 2$.*

PROOF. From the shape of the probability generating functional of a Cox process directed by an arbitrary random measure [see Daley and Vere-Jones (1988), page 262], we deduce that

$$G_{\Pi_\Xi^{(z)}}[w] = \mathsf{E}\exp\left\{z \int_{R^d} (w(x) - 1)\mathbf{1}_{\Xi^c}(x)\, dx\right\},$$

which in turn, using the above definitions of moment and cumulant measures, provides

$$\alpha_{\Pi_\Xi^{(z)}}^{(n)}\left(\underset{j=1}{\overset{n}{\times}} B_j\right) = z^n \mathsf{E} \prod_{j=1}^n |\Xi^c \cap B_j|$$

and

$$\gamma_{\Pi_\Xi^{(z)}}^{(n)}\left(\underset{j=1}{\overset{n}{\times}} B_j\right) = z^n \Gamma(|\Xi^c \cap B_1|, \ldots, |\Xi^c \cap B_n|).$$

Hence, repeating the steps in the proof of Lemma 1 that lead to (2.3) and (2.5) (with $\Xi^c$ instead of $\Xi$), we recognize that, for $B \in \mathfrak{B}(R^{dn})$,

$$\alpha^{(n)}_{\Pi^{(z)}_\Xi}(B) = z^n \int_B p^{(n)}_{\Xi^c}(X_n)\,dX_n \quad \text{and} \quad \gamma^{(n)}_{\Pi^{(z)}_\Xi}(B) = z^n \int_B c^{(n)}_{\Xi^c}(X_n)\,dX_n,$$

where $X_n = \{x_1, \ldots, x_n\}$ and $dX_n = d(x_1, \ldots, x_n)$. Thus, the $n$th-order cumulant density of $\Pi^{(z)}_\Xi$ equals $z^n c^{(n)}_{\Xi^c}$. The proof is completed by appealing to (2.6), (2.7) and the very definition of total variation. □

We now introduce a further family of (signed) measures $\gamma^{(m,n)}_\Psi$ on $\mathfrak{B}(R^{m+n})$ for $n, m \geq 0$ associated with the point process $\Psi$, which is assumed to admit moment measures of order $m + n$. For bounded $A_1, \ldots, A_m, B_1, \ldots, B_n \in \mathfrak{B}(R^d)$ define

$$\gamma^{(m,n)}_\Psi\left(\underset{i=1}{\overset{m}{\times}} A_i \times \underset{j=1}{\overset{n}{\times}} B_j\right)$$
$$:= \lim_{\substack{u_1,\ldots,u_m\uparrow 1 \\ v_1,\ldots,v_n\uparrow 1}} \frac{\partial^{n+m}}{\partial v_1 \cdots \partial v_n\, \partial u_1 \cdots \partial u_m} \frac{G_\Psi[w^{A_1,\ldots,A_m,B_1,\ldots,B_n}_{u_1,\ldots,u_m,v_1,\ldots,v_n}]}{G_\Psi[w^{B_1,\ldots,B_n}_{v_1,\ldots,v_n}]}.$$

For the sake of distinction, let us call $\gamma^{(m,n)}_\Psi$ the (factorial) *correlation measure of order* $(m, n)$. In case the moment density $p^{(m+n)}_\Psi$ exists, $\gamma^{(m,n)}_\Psi$ has a Lebesgue density $c^{(m,n)}_\Psi$ which we call *correlation density of order* $(m, n)$. Note that Ruelle (1964) apparently first introduced the densities $c^{(m,n)}_\Psi$ via an algebraic method to study cluster properties of the correlation functions of classical gases.

It is evident that $\gamma^{(m,n)}_\Psi$ is symmetric in the first $m$ as well as in the second $n$ components, but not completely symmetric. By logarithmic differentiation with respect to $v_1$ in the above definition of $\gamma^{(n)}_\Psi$ we see that $\gamma^{(n)}_\Psi$ coincides with $\gamma^{(1,n-1)}_\Psi$ for $n \geq 1$ and, thus, $c^{(n)}_\Psi = c^{(1,n-1)}_\Psi$ for $n \geq 1$ provided the densities exist. Moreover, for fixed $m \geq 1$ and any $n \geq 1$, the relationship between factorial moment and correlation measures,

$$(2.8) \quad \alpha^{(m+n)}_\Psi\left(\underset{i=1}{\overset{m}{\times}} A_i \times \underset{j=1}{\overset{n}{\times}} B_j\right)$$
$$= \sum_{\substack{\varnothing \subseteq J \subseteq N \\ N=\{1,\ldots,n\}}} \gamma^{(m,|J|)}_\Psi\left(\underset{i=1}{\overset{m}{\times}} A_i \times \underset{j \in J}{\times} B_j\right) \alpha^{(n-|J|)}_\Psi\left(\underset{j \in N\setminus J}{\times} B_j\right),$$

holds, where the summation extends over all subsets $J$ of $N$ with $|J|$ elements. For reasons of consistency, put $\gamma^{(m,0)}_\Psi(\times^m_{i=1} A_i) = \alpha^{(m)}_\Psi(\times^m_{i=1} A_i)$ ($= 1$ for $m = 0$) and $\gamma^{(0,n)}_\Psi(\times^n_{j=1} B_j) = 0$ for $n \geq 1$.

To verify (2.8), let us briefly write $D_v$ for $\partial^n/\partial v_1 \cdots \partial v_n$ and $D_u$ for $\partial^m/\partial u_1 \cdots \partial u_n$ and put

$$f(u,v) = w^{A_1,\ldots,A_m,B_1,\ldots,B_n}_{u_1,\ldots,u_m,v_1,\ldots,v_n} \quad \text{and} \quad g(v) = w^{B_1,\ldots,B_n}_{v_1,\ldots,v_n}.$$

Then (2.8) is obtained by applying Leibniz's rule for higher-order derivatives of products of functions to the right-hand side of the identity

$$D_v D_u(f(u,v)) = D_v(g(v) \cdot h(v)) \qquad \text{with } h(v) = \frac{D_u(f(u,v))}{g(v)}.$$

We conclude this section with a recursive representation of the correlation density $c_\Psi^{(m,n)}$ in terms of the densities $c_\Psi^{(m-1+j,n-j)}$, $j = 1,\ldots,n$. For notational convenience, we omit the superscripts (if confusion is excluded) and write $c_\Psi(X_m, Y_n)$ instead of $c_\Psi^{(m,n)}(x_1,\ldots,x_m,y_1,\ldots,y_n)$, where $X_m = \{x_1,\ldots,x_m\}$ and $Y_n = \{y_1,\ldots,y_n\}$ are two disjoint sets of distinct points in $R^d$. Furthermore, put $X'_{m-1} = X_m \setminus \{x_1\}$ and let $|Y|$ denote the cardinality of a finite point set $Y \subset R^d$.

LEMMA 3. *Let $\Psi$ be a point process on $R^d$ with strictly positive factorial moment densities $p_\Psi^{(k)}$ for $k = 1,\ldots,m+n \geq 1$. Then we have*

$$(2.9) \quad c_\Psi(X_m, Y_n) = \sum_{\varnothing \subseteq Y \subseteq Y_n} (-1)^{|Y|} K_\Psi(X_m, Y) c_\Psi(Y \cup X'_{m-1}, Y_n \setminus Y),$$

*where*

$$K_\Psi(X_m, Y) := \sum_{\varnothing \subseteq V \subseteq Y} (-1)^{|V|} \frac{p_\Psi(V \cup X_m)}{p_\Psi(V \cup X'_{m-1})}$$

(2.10)
$$\text{for } m, |Y| \geq 1, Y \subseteq Y_n$$

*and $K_\Psi(X_m, \varnothing) = p_\Psi(X_m)/p_\Psi(X'_{m-1})$ for $m \geq 1$ and $K_\Psi(\varnothing, Y_n) = 0$ for $n \geq 1$.*

PROOF. The relationship (2.8) reads, in terms of densities, as

$$(2.11) \qquad p_\Psi(X_m \cup Y_n) = \sum_{\varnothing \subseteq Y \subseteq Y_n} c_\Psi(X_m, Y) p_\Psi(Y_n \setminus Y).$$

Given the moment density functions $p_\Psi(Y)$, $Y \subseteq Y_n$, with $p_\Psi(\varnothing) = 1$, there exist unique symmetric functions $p_\Psi^*(Y)$, $Y \subseteq Y_n$, with $p_\Psi^*(\varnothing) = 1$ that satisfy the equations

$$(2.12) \qquad \sum_{\varnothing \subseteq V \subseteq Y} p_\Psi^*(V) p_\Psi(Y \setminus V) = 0 \qquad \text{for } \varnothing \neq Y \subseteq Y_n.$$

By means of the functions $p_\Psi^*(Y)$ we may invert the "convolution equation" (2.11) by calculating the sum

$$\sum_{\varnothing \subseteq Y \subseteq Y_n} p_\Psi^*(Y_n \setminus Y) p_\Psi(X_m \cup Y)$$
$$= \sum_{\varnothing \subseteq Y \subseteq Y_n} p_\Psi^*(Y_n \setminus Y) \sum_{\varnothing \subseteq V \subseteq Y} c_\Psi(X_m, V) p_\Psi(Y \setminus V)$$
$$= \sum_{\varnothing \subseteq V \subseteq Y_n} c_\Psi(X_m, V) \sum_{Y:\, V \subseteq Y \subseteq Y_n} p_\Psi^*(Y_n \setminus Y) p_\Psi(Y \setminus V).$$

Since, by (2.12), the second sum in the last line vanishes for all proper subsets $V \subset Y_n$, the whole last line is equal to $c_\Psi(X_m, Y_n)$. Using this identity and the relationship

$$\sum_{\varnothing \subseteq V \subseteq Y} (-1)^{|V|} K_\Psi(X_m, V) = \frac{p_\Psi(Y \cup X_m)}{p_\Psi(Y \cup X'_{m-1})} \qquad \text{for } \varnothing \subseteq Y \subseteq Y_n$$

obtained from (2.10) by using the Möbius inversion formula [see Rota (1964)], we may proceed with

$$c_\Psi(X_m, Y_n)$$
$$= \sum_{\varnothing \subseteq Y \subseteq Y_n} p_\Psi^*(Y_n \setminus Y) p_\Psi(Y \cup X'_{m-1}) \sum_{\varnothing \subseteq V \subseteq Y} (-1)^{|V|} K_\Psi(X_m, V)$$
$$= \sum_{\varnothing \subseteq V \subseteq Y_n} (-1)^{|V|} K_\Psi(X_m, V) \sum_{Y:\, V \subseteq Y \subseteq Y_n} p_\Psi^*(Y_n \setminus Y) p_\Psi(Y \cup X'_{m-1})$$
$$= \sum_{\varnothing \subseteq V \subseteq Y_n} (-1)^{|V|} K_\Psi(X_m, V)$$
$$\quad \times \sum_{\varnothing \subseteq U \subseteq Y_n \setminus V} p_\Psi^*(Y_n \setminus V \setminus U) p_\Psi(U \cup V \cup X'_{m-1}).$$

Applying again the above derived identity, we see that the second sum in the last line equals $c_\Psi(V \cup X'_{m-1}, Y_n \setminus V)$, proving the asserted relationship (2.9). □

**3. Absolute integrability of the correlation densities of the Cox process $\Pi_\Xi^{(1)}$ and proof of Theorem 1.** Throughout this section we consider the factorial moment and correlation densities $p_\Psi$ and $c_\Psi$ merely with respect to Cox process $\Pi_\Xi^{(1)}$ defined by (2.1) for the PGM (1.1). For notational ease, we indicate this by omitting the subscript $\Psi$ at $p_\Psi$, $c_\Psi$ and $K_\Psi$. Our aim is to obtain bounds of the integrals $\int_{R^{dn}} |c(\{o\}, Y_n)|\, dY_n$ [$= G_n(\Xi)$ by Lemma 2] under suitable moment conditions on $|\Xi_0|$. To do this, however, our inductive proving technique requires us to estimate the integrals $\int_{R^{dn}} |c(X_m, Y_n)| \, dY_n$ uniformly in $X_m \in R^{dm}$ for any $m \geq 1$.

Let $X_m$, $X'_{m-1}$ and $Y_n$ be the finite point sets introduced at the end of Section 2. Furthermore, for any finite subset $Y \subset R^d$, put $\Xi_0(Y) := \bigcup_{y \in Y}(\Xi_0 - y)$ [and thus $\Xi_0^c(Y) = \bigcap_{y \in Y}(\Xi_0^c - y)$]. For any $n \geq 1$ and $Y \subseteq Y_n$, define

$$(3.1) \qquad S(X_m, Y) := \sum_{\varnothing \subseteq V \subseteq Y} (-1)^{|V|} \exp\{E(x_1; X'_{m-1}, V)\},$$

where

$$E(x; U, V) := \lambda \mathsf{E}|(\Xi_0 - x) \cap \Xi_0^c(U) \cap \Xi_0(V)| \qquad \text{for } x \notin U$$

and, for any $Y \subseteq Y_{n-1} := Y_n \setminus \{y_n\}$,

$$(3.2) \qquad T(y_n, X_m, Y) := \sum_{\varnothing \subseteq V \subseteq Y} (-1)^{|V|} \exp\{-E(x_1, y_n; X'_{m-1}, V)\},$$

where $E(x, y; U, V) := \lambda \mathsf{E}|(\Xi_0 - x) \cap (\Xi_0 - y) \cap \Xi_0^c(U) \cap \Xi_0(V)|$ for $x \neq y$, $x \notin U$ and $y \notin V$. Note that $E(x, y; U, \varnothing) = E(x; U, \varnothing) = 0$, implying $S(X_m, \varnothing) = T(y_n, X_m, \varnothing) = 1$.

From (1.6) it is clear that

$$\frac{p(V \cup X_m)}{p(V \cup X'_{m-1})} = \exp\{-\lambda \mathsf{E}|\Xi_0| + \lambda \mathsf{E}|(\Xi_0 - x_1) \cap \Xi_0(V \cup X'_{m-1})|\}$$

$$= \exp\{-\lambda \mathsf{E}|(\Xi_0 - x_1) \cap \Xi_0^c(X'_{m-1})|\} \exp\{E(x_1; X'_{m-1}, V)\},$$

so that, by (2.10) and (3.1),

$$(3.3) \qquad K(X_m, Y) = \exp\{-\lambda \mathsf{E}|(\Xi_0 - x_1) \cap \Xi_0^c(X'_{m-1})|\} S(X_m, Y).$$

Next we establish a recursive representation of $S(X_m, Y_n)$ with respect to $Y_n$ in combination with the nonnegative terms $T(y_n, X_m, Y)$ for $Y \subseteq Y_{n-1}$. It turns out (see Lemma 5 below) that the integrals $\int_{R^{dn}} T(y_n, X_m, Y_{n-1}) \, dY_n$ can be represented as functionals of certain PGM (3.6), which enables us to derive upper bounds of them under reasonable moment conditions on the volume of the typical grain $\Xi_0$. By means of these bounds and the following Lemma 4 we find corresponding bounds of $\int_{R^{dn}} |S(X_m, Y_n)| \, dY_n$ which in turn, using (2.9) with (3.3), enable us to establish the desired bounds of $\int_{R^{dn}} |c(X_m, Y_n)| \, dY_n$; see Lemma 7 below.

LEMMA 4. *We have*

$$S(X_m, Y_n) = S(X_m, Y_{n-1})(1 - \exp\{E(x_1; X'_{m-1}, \{y_n\})\})$$

$$- \exp\{E(x_1; X'_{m-1}, \{y_n\})\}$$

$$\times \sum_{\varnothing \subset Y \subseteq Y_{n-1}} T(y_n, X_m, Y)$$

$$\times \exp\{E(x_1; X'_{m-1}, Y)\} S(X_m \cup Y, Y_{n-1} \setminus Y).$$

PROOF. By the definition of the terms $E(x;U,V)$ and $E(x,y;U,V)$, and the relationship $|A|+|B|-|A\cap B|=|A\cup B|$ for bounded $A,B\in\mathfrak{B}(R^d)$, we get

$$E(x_1;X'_{m-1},Y\cup\{y_n\})$$
$$= E(x_1;X'_{m-1},\{y_n\})+E(x_1;X'_{m-1},Y)-E(x_1,y_n;X'_{m-1},Y)$$

for any $Y\subseteq Y_{n-1}$. Furthermore, we may rewrite the sum $S(X_m,Y_n)$ as

$$\sum_{\varnothing\subseteq Y\subseteq Y_{n-1}}(-1)^{|Y|}(\exp\{E(x_1;X'_{m-1},Y)\}-\exp\{E(x_1;X'_{m-1},Y\cup\{y_n\})\}).$$

This combined with the foregoing relationship leads to

$$S(X_m,Y_n) = S(X_m,Y_{n-1})$$
$$- \exp\{E(x_1;X'_{m-1},\{y_n\})\}$$
$$\times \sum_{\varnothing\subseteq Y\subseteq Y_{n-1}}(-1)^{|Y|}\exp\{E(x_1;X'_{m-1},Y)-E(x_1,y_n;X'_{m-1},Y)\}.$$

A simple application of the Möbius inversion formula [see Rota (1964)] to the terms (3.2) yields

$$\exp\{-E(x_1,y_n;X'_{m-1},Y)\}=1+\sum_{\varnothing\subset U\subseteq Y}(-1)^{|U|}T(y_n,X_m,U),\qquad Y\subseteq Y_{n-1}.$$

Inserting this identity on the right-hand side of the previous equality we arrive at

$$S(X_m,Y_n)$$
$$= S(X_m,Y_{n-1})(1-\exp\{E(x_1;X'_{m-1},\{y_n\})\})$$
$$- \exp\{E(x_1;X'_{m-1},\{y_n\})\}$$
$$\times \sum_{\varnothing\subset Y\subseteq Y_{n-1}}\sum_{\varnothing\subset U\subseteq Y}(-1)^{|Y|-|U|}T(y_n,X_m,U)\exp\{E(x_1;X'_{m-1},Y)\}.$$

By interchanging the sums and substituting $V=Y\setminus U$ we obtain that

$$\sum_{\varnothing\subset Y\subseteq Y_{n-1}}\sum_{\varnothing\subset U\subseteq Y}(-1)^{|Y|-|U|}T(y_n,X_m,U)\exp\{E(x_1;X'_{m-1},Y)\}$$
$$= \sum_{\varnothing\subset U\subseteq Y_{n-1}}\sum_{Y:U\subseteq Y\subseteq Y_{n-1}}(-1)^{|Y|-|U|}T(y_n,X_m,U)$$
$$\times \exp\{E(x_1;X'_{m-1},(Y\setminus U)\cup U\}$$
$$= \sum_{\varnothing\subset U\subseteq Y_{n-1}}T(y_n,X_m,U)\sum_{\varnothing\subseteq V\subseteq Y_{n-1}\setminus U}(-1)^{|V|}\exp\{E(x_1;X'_{m-1},V\cup U)\}.$$

Since $|A \cup B| = |B| + |A \cap B^c|$ for any bounded $A, B \in \mathfrak{B}(R^d)$,

$$E(x_1; X'_{m-1}, V \cup U) = E(x_1; X'_{m-1}, U) + E(x_1; U \cup X'_{m-1}, V),$$

whence, by definition (3.1), it follows that

$$\sum_{\varnothing \subseteq V \subseteq Y_{n-1} \setminus U} (-1)^{|V|} \exp\{E(x_1; X'_{m-1}, V \cup U)\}$$
$$= \exp\{E(x_1; X'_{m-1}, U)\} S(X_m \cup U, Y_{n-1} \setminus U).$$

Finally, assembling all the above identities we obtain the assertion of Lemma 4. □

LEMMA 5. *Let $\Xi$ be the PGM (1.1) with compact typical grain $\Xi_0$ satisfying $\mathsf{E}|\Xi_0|^{n+1} < \infty$ for some fixed $n \geq 2$. Then, for any $m \geq 1$,*

$$\sup_{X_m} \int_{R^{dn}} T(y_n, X_m, Y_{n-1}) \, dY_n$$

(3.4)
$$\leq (n-1)!$$
$$\times \sum_{k=1}^{n-1} \frac{\lambda^k}{k!} \sum_{\substack{n_1+n_2+\cdots+n_k=n-1 \\ n_1,\ldots,n_k \geq 1}} \frac{\mathsf{E}|\Xi_0|^{n_1+2}}{n_1!} \frac{\mathsf{E}|\Xi_0|^{n_2+1}}{n_2!} \cdots \frac{\mathsf{E}|\Xi_0|^{n_k+1}}{n_k!}.$$

*If condition (1.3) is satisfied, then the estimate*

$$(3.5) \quad \sup_{X_m} \int_{R^{dn}} T(y_n, X_m, Y_{n-1}) \, dY_n \leq n! \left(\frac{2}{a}\right)^n \frac{\lambda M(a)}{a} \left(1 + \frac{\lambda M(a)}{a}\right)^{n-2}$$

*holds for all $n \geq 2$ and $m \geq 1$.*

PROOF. According to the definition (3.2),

$$T(y_n, X_m, Y_{n-1}) = \sum_{\varnothing \subseteq Y \subseteq Y_{n-1}} (-1)^{|Y|} \exp\{-E(x_1, y_n; X'_{m-1}, Y)\}.$$

Now, for any nonempty $Y \subseteq Y_{n-1}$, we introduce a new PGM $\Xi(x_1, y_n; X'_{m-1}, Y)$ governed by $\Pi_\lambda = \sum_{i \geq 1} \delta_{X_i}$ and the typical grain $(\Xi_0 - x_1) \cap (\Xi_0 - y_n) \cap \Xi_0^c(X'_{m-1}) \cap \Xi_0(Y)$, that is,

(3.6)
$$\Xi(x_1, y_n; X'_{m-1}, Y)$$
$$:= \bigcup_{i \geq 1} ((\Xi_i - x_1) \cap (\Xi_i - y_n) \cap \Xi_i^c(X'_{m-1}) \cap \Xi_i(Y) + X_i),$$

where $\Xi_i(Y) = \bigcup_{y \in Y} (\Xi_i - y)$ and $\Xi_i^c(X'_{m-1}) = \bigcap_{j=2}^m (\Xi_i^c - x_j)$.

Obviously, for each realization of the PGM (3.6), we have

$$\Xi(x_1, y_n; X'_{m-1}, Y) = \bigcup_{y \in Y} \Xi(x_1, y_n; X'_{m-1}, \{y\}) \qquad \text{for } Y \subseteq Y_{n-1}.$$

Applying the well-known formula $\mathsf{P}(o \in \Xi) = 1 - \exp\{-\lambda \mathsf{E}|\Xi_0|\}$ [which is valid for the PGM (1.1)] to the stationary PGM (3.6) we see that

$$\exp\{-E(x_1, y_n; X'_{m-1}, Y)\} = \mathsf{P}(o \notin \Xi(x_1, y_n; X'_{m-1}, Y))$$

$$= 1 - \mathsf{P}\Biggl(\bigcup_{y \in Y} \{o \in \Xi(x_1, y_n; X'_{m-1}, \{y\})\}\Biggr).$$

Since $\sum_{\varnothing \subseteq Y \subseteq Y_{n-1}} (-1)^{|Y|} = 0$, it follows from the inclusion–exclusion principle that

$$T(y_n, X_m, Y_{n-1}) = - \sum_{\varnothing \subset Y \subseteq Y_{n-1}} (-1)^{|Y|} \mathsf{P}\Biggl(\bigcup_{y \in Y} \{o \in \Xi(x_1, y_n; X'_{m-1}, \{y\})\}\Biggr)$$

$$= \mathsf{P}\Biggl(\bigcap_{i=1}^{n-1} \{o \in \Xi(x_1, y_n; X'_{m-1}, \{y_i\})\}\Biggr)$$

$$= \mathsf{E}\Biggl(\prod_{i=1}^{n-1} \mathbf{1}_{\Xi(x_1, y_n; X'_{m-1}, \{y_i\})}(o)\Biggr).$$

Thus, by Fubini's theorem,

$$\int_{R^{dn}} T(y_n, X_m, Y_{n-1}) \, dY_n = \int_{R^d} \mathsf{E}\Biggl(\int_{R^d} \mathbf{1}_{\Xi(x_1, y_n; X'_{m-1}, \{y\})}(o) \, dy\Biggr)^{n-1} dy_n$$

and, for each realization of (3.6),

$$\int_{R^d} \mathbf{1}_{\Xi(\{x_1, y_n\}, X'_{m-1}, \{y\})}(o) \, dy \leq \int_{R^d} \sum_{i \geq 1} \mathbf{1}_{(\Xi_i - x_1) \cap (\Xi_i - y_n) \cap (\Xi_i - y)}(X_i) \, dy$$

$$\leq \sum_{i \geq 1} |\Xi_i| \mathbf{1}_{(\Xi_i - x_1) \cap (\Xi_i - y_n)}(X_i),$$

whence, by applying the polynomial formula and using Fubini's theorem again, we get that

$$\mathsf{E}\Biggl(\int_{R^d} \mathbf{1}_{\Xi(x_1, y_n; X'_{m-1}, \{y\})}(o) \, dy\Biggr)^{n-1}$$

$$\leq \sum_{k=1}^{n-1} \frac{1}{k!} \sum_{\substack{n_1 + \cdots + n_k = n-1 \\ n_1, \ldots, n_k \geq 1}} \frac{(n-1)!}{n_1! \cdots n_k!}$$

$$\times \mathsf{E}\bigg(\sum_{i_1,\ldots,i_k\geq 1}^{*}\prod_{j=1}^{k}(\mathbf{1}_{(\Xi_{i_j}-x_1)\cap(\Xi_{i_j}-y_n)}(X_{i_j})|\Xi_{i_j}|^{n_j})\bigg)$$

$$=(n-1)!\sum_{k=1}^{n-1}\frac{\lambda^k}{k!}\sum_{\substack{n_1+\cdots+n_k=n-1\\n_1,\ldots,n_k\geq 1}}\prod_{j=1}^{k}\frac{\mathsf{E}|(\Xi_0-x_1)\cap(\Xi_0-y_n)||\Xi_0|^{n_j}}{n_j!}.$$

Here the sum $\sum^*$ stretches over $k$-tuples of pairwise distinct indices and the last equality is obtained by applying the Campbell-type formula (A.3) for $f_j(x,K)=\mathbf{1}_{(K-x_1)\cap(K-y_n)}(x)|K|^{n_j}$. Together with the obvious relationship

$$(3.7)\qquad \int_{R^d}|(\Xi_0-x_1)\cap(\Xi_0-y_n)|\,dy_n=|\Xi_0|^2,$$

we finally arrive at the desired estimate (3.4).

The existence of the exponential moment $M(a)$ of $|\Xi_0|$ implies $\mathsf{E}|\Xi_0|^k\leq k!M(a)a^{-k}$ for all $k\geq 1$. Inserting this moment bound in the right-hand side of (3.4) and taking into account

$$\sum_{\substack{n_1+\cdots+n_k=n-1\\n_1,\ldots,n_k\geq 1}}(n_1+2)\prod_{i=1}^{k}(n_i+1)\leq\binom{n-2}{k-1}2^n n,$$

we obtain that

$$\int_{R^d}\mathsf{E}\bigg(\int_{R^d}\mathbf{1}_{\Xi(x_1,y_n;X_m\setminus\{x_1\},\{y\})}(o)\,dy\bigg)^{n-1}dy_n$$

$$\leq\frac{n!}{a^n}\sum_{k=1}^{n-1}\frac{(\lambda M(a))^k}{a^k k!}\binom{n-2}{k-1}2^n$$

$$\leq n!\bigg(\frac{2}{a}\bigg)^n\frac{\lambda M(a)}{a}\bigg(1+\frac{\lambda M(a)}{a}\bigg)^{n-2}.$$

This is exactly the desired estimate (3.5). Thus, Lemma 5 is completely proved.

$\square$

LEMMA 6. *Let $\Xi$ be the PGM (1.1) with compact typical grain $\Xi_0$ that satisfies $\mathsf{E}|\Xi_0|^{n+1}<\infty$ for some fixed $n\geq 1$. Then, for any $m\geq 1$,*

$$(3.8)\qquad \sup_{X_m}\int_{R^{dn}}|S(X_m,Y_n)|\,dY_n\leq c_n(\lambda)<\infty,$$

*where the constant $c_n(\lambda)$ depends on $\lambda$ and the first $n+1$ moments of $|\Xi_0|$. Moreover, if condition (1.3) is satisfied, then (3.8) holds with*

$$(3.9)\qquad c_n(\lambda)=n!A^n B(1+B)^{n-1}$$

for all $n \geq 1$ and $m \geq 1$ with $A = \frac{2}{a}(1 + \exp\{\lambda\mathsf{E}|\Xi_0|\}) = \frac{2(2-p)}{a(1-p)}$ and $B = \frac{\lambda M(a)}{a}$.

PROOF. In view of the obvious inequalities
$$E(x_1; X'_{m-1}, \{y_1\}) \leq \mathsf{E}|(\Xi_0 - x_1) \cap (\Xi_0 - y_1)| \leq \lambda\mathsf{E}|\Xi_0|$$
and $e^x - 1 \leq xe^x$ for $x \geq 0$ together with (3.7) we see that
$$\int_{R^d} |S(X_m, \{y_1\})| \, dy_1 = \int_{R^d} (\exp\{E(x_1; X'_{m-1}, \{y_1\})\} - 1) \, dy_1$$
$$\leq \lambda \exp\{\lambda\mathsf{E}|\Xi_0|\} \int_{R^d} \mathsf{E}|(\Xi_0 - x_1) \cap (\Xi_0 - y_1)| \, dy_1$$
$$= \frac{\lambda}{1-p} \mathsf{E}|\Xi_0|^2.$$

Define
$$A_{m,n} := \sup_{X_m} \int_{R^{dn}} |S(X_m, Y_n)| \, dY_n \quad \text{and} \quad B_{m,n} := \sup_{X_m} \int_{R^{dn}} T(y_n, X_m, Y_{n-1}) \, dY_n$$

for $m, n \geq 1$. From Lemma 4 and $E(x_1; X'_{m-1}, Y) \leq \lambda\mathsf{E}|\Xi_0|$ for $Y \subseteq Y_{n-1}$ we get the inequality

(3.10)
$$A_{m,n} \leq A_{m,n-1} \frac{\lambda}{1-p} \mathsf{E}|\Xi_0|^2$$
$$+ \frac{1}{(1-p)^2} \sum_{k=1}^{n-1} \binom{n-1}{k} B_{m,k+1} A_{m+k,n-k-1}$$

with $A_{m,0} = 1$ and $A_{m,1} \leq \lambda\mathsf{E}|\Xi_0|^2/(1-p)$ for any $m \geq 1$. Since, by Lemma 5, $B_{m,k+1} \leq k!C_k$ for $k \geq 1$, where $C_k$ depends on $\lambda$ and the first $k+2$ moments of $|\Xi_0|$ but not on $m$, we recognize by induction on $n$, that $A_{m,n} \leq n!D_n$, where $D_0 = 1$, $D_1 = \lambda\mathsf{E}|\Xi_0|^2/(1-p)$ and
$$D_n = \frac{1}{(1-p)n}\left(D_{n-1}C_0 + \frac{1}{1-p}\sum_{k=1}^{n-1} C_k D_{n-k-1}\right)$$

for $n \geq 2$ with $C_0 := \lambda\mathsf{E}|\Xi_0|^2$. Therefore, $A_{m,n}$ does not depend on $m$ and is bounded by terms that involve merely $\lambda$ and $\mathsf{E}|\Xi_0|^k$, $k = 1, 2, \ldots, n+1$. This proves the first part of Lemma 6.

We also prove (3.9) by induction on $n$. From (1.3) we get $\mathsf{E}|\Xi_0|^2 \leq 2M(a)/a^2$, implying
$$A_{m,1} \leq D_1 \leq \lambda\exp\{\lambda\mathsf{E}|\Xi_0|\}\mathsf{E}|\Xi_0|^2 \leq \frac{2\lambda M(a)}{a^2(1-p)},$$

which is even slightly stronger than (3.9) for $n=1$. Assume now the validity of (3.9) for $n=1,\ldots,N-1$. Taking into account the estimates $B_{m,k+1} \leq k!C_k$ with $C_k = (k+1)B(\frac{2}{a})^{k+1}(1+B)^{k-1}$ for $k \geq 1$ as stated in Lemma 5 together with $C_0 \leq 2B/a$, we may write

$$D_N \leq \frac{1}{N(1-p)}$$
$$\times \left( A^{N-1}B(1+B)^{N-2}\frac{2B}{a} \right.$$
$$\left. + \frac{1}{1-p}\sum_{k=1}^{N-1}(k+1)B\left(\frac{2}{a}\right)^{k+1}(1+B)^{k-1}A^{N-k-1}(1+B)^{N-k-1} \right).$$

After a short calculation using that $\sum_{k\geq 1}(k+1)(\frac{2}{aA})^k = (1-p)(3-p)$, we arrive at

$$D_N \leq A^N B(1+B)^{N-1} \qquad \text{for } N \geq 2.$$

Thus, the second part of Lemma 6 is proved. $\square$

LEMMA 7. *Let $\Xi$ be the PGM (1.1) with compact typical grain $\Xi_0$ that satisfies $\mathsf{E}|\Xi_0|^{n+1} < \infty$ for some fixed $n \geq 1$. Then, for any $m \geq 1$,*

$$\sup_{X_m} \int_{R^{dn}} |c(X_m, Y_n)|\, dY_n \leq c_{m,n}(\lambda) < \infty, \tag{3.11}$$

*where the constant $c_{m,n}(\lambda)$ depends on $m$, $\lambda$ and the first $n+1$ moments of $|\Xi_0|$. If condition (1.3) is satisfied, then (3.11) holds with*

$$c_{m,n}(\lambda) = n!\, 2^{m+1} AB(4A(1+B))^{n-1} \tag{3.12}$$

*for all $n \geq 1$ and $m \geq 1$ with $A$ and $B$ as in Lemma 6.*

PROOF. Replacing $K(X_m, Y)$ in (2.9) with (1.6) leads to

$$c(X_m, Y_n) = \exp\{-\lambda \mathsf{E}|(\Xi_0 - x_1) \cap \Xi_0^c(X'_{m-1})|\}$$
$$\times \sum_{\varnothing \subseteq Y \subseteq Y_n} (-1)^{|Y|} S(X_m, Y) c(Y \cup X'_{m-1}, Y_n \setminus Y).$$

Since $c(X_m, \varnothing) = \exp\{-\lambda|\Xi_0(X_m)|\} \leq 1$ by (1.6), $c(\varnothing, Y_n) = 0$ for $n \geq 1$ by definition and $S(X_m, \varnothing) = 1$ for $m \geq 1$, and since both $S(X,Y)$ and $c(X,Y)$ are symmetric in $Y \subseteq Y_n$ for fixed $X$, we deduce from the latter recurrence relationship the inequality

$$\int_{R^{dn}} |c(X_m, Y_n)|\, dY_n$$

$$
\begin{aligned}
&\leq \int_{R^{dn}} |c(X'_{m-1}, Y_n)|\, dY_n + \int_{R^{dn}} |S(X_m, Y_n)|\, dY_n \\
&\quad + \sum_{k=1}^{n-1} \binom{n}{k} \int_{R^{dk}} |S(X_m, Y_k)|\, dY_k \\
&\qquad \times \sup_{Y_k} \int_{R^{d(n-k)}} |c(Y_k \cup X'_{m-1}, Y_n \setminus Y_k)|\, d(Y_n \setminus Y_k).
\end{aligned}
\tag{3.13}
$$

For any $m \geq 1$ we have

$$
\begin{aligned}
\int_{R^d} |c(X_m, \{y\})|\, dy &= \int_{R^d} \exp(-\lambda \mathsf{E}|\Xi_0(X_m) \cup (\Xi_0 - y)|) \\
&\quad \times (1 - \exp(-\lambda \mathsf{E}|\Xi_0(X_m) \cap (\Xi_0 - y)|))\, dy \\
&\leq \lambda(1-p)\mathsf{E}|\Xi_0(X_m)||\Xi_0| \leq \lambda m \mathsf{E}|\Xi_0|^2.
\end{aligned}
$$

Using the estimate (3.8) of Lemma 6 and applying (3.13) successively to the remaining integrals on the right-hand side of (3.13), we obtain a bound of the left-hand side of (3.13) in terms of $c_k(\lambda)$, $k = 1, \ldots, n$, and $\sup_Y \int_{R^d} |c(X \cup Y, \{y_n\})|\, dy_n$, $Y \subseteq Y_{n-1}$, $X \subseteq X'_{m-1}$. This combined with the foregoing inequality proves (3.11).

We now assume (1.3), which gives $\mathsf{E}|\Xi_0|^2 \leq 2M(a)/a^2$, so that together with $m \leq 2^{m-1}$,

$$
\int_{R^d} |c(X_m, \{y\})|\, dy \leq 2^{m-1} AB,
$$

which implies (3.12) for $n = 1$ and $m \geq 1$. Let now (3.12) hold for all $m, n \geq 1$ that satisfy $m + n < M + N$. Then, making use of estimate (3.9) of Lemma 6, it follows from (3.13) that

$$
\begin{aligned}
&\int_{R^{dN}} |c(X_M, Y_N)|\, dY_N \\
&\leq N! 2^{M-2} (4A)^N B(1+B)^{N-1} + N! A^N B(1+B)^{N-1} \\
&\quad + N! \sum_{k=1}^{N-1} A^k B(1+B)^{k-1} 2^{M+k-2} (4A(1+B))^{N-k} \\
&= N! 2^{M-1} (4A)^N B(1+B)^{N-1} \left[ \frac{1}{2} + \frac{1}{2^{2N+M-1}} + \sum_{k=1}^{N-1} \frac{1}{2^{k+1}} \right].
\end{aligned}
$$

Thus, the validity of (3.12) for $m + n = M + N$ follows because the sum in brackets does not exceed one for $M + N \geq 2$. This completes the proof of Lemma 7. $\square$

PROOF OF THEOREM 1. As an immediate consequence of (3.12) for $m = 1$ and Lemma 2 applied to the stationary PGM (1.1), we obtain that

$$(3.14) \quad G_{n+1}(\Xi) = \int_{R^{dn}} |c(\{o\}, Y_n)| \, dY_n \leq c_{1,n}(\lambda) \leq n! 4AB(4A(1+B))^{n-1}$$

for all $n \geq 1$. Thus, by the definition of $A$ and $B$ in Lemma 5, we get (1.10) with $H(a) = 4AB$ and $\Delta(a) = 4A(1 + B)$. Finally, the existence and analyticity of the thermodynamic limit $L(z)$ of the function (1.2) on the disk $D_{\Delta(a)}$ follows from the second part of Lemma 1. □

## 4. Proofs of Theorems 2–4.

PROOF OF THEOREM 2. From Lemma 1 and (3.14) we get the estimate

$$(4.1) \qquad |\Gamma_k(|\Xi \cap W|)| \leq |W|(k-1)! H(a) \Delta(a)^{k-2} \qquad \text{for } k \geq 2$$

and any $W \in \mathfrak{B}(R^d)$. For the standardized random variable $\xi_n := (|\Xi \cap W_n| - p|W_n|)/\sigma_n \sqrt{|W_n|}$ (with $\sigma_n > 0$), (4.1) implies that

$$(4.2) \qquad |\Gamma_k(\xi_n)| \leq (k-1)! \frac{H(a)}{\sigma_n^2 \Delta_n^{k-2}} \leq \frac{k! H_n}{\Delta_n^{k-2}} \qquad \text{for } k \geq 3$$

with $H_n = H(a)/2\sigma_n^2$ and $\Delta_n = \sigma_n \sqrt{|W_n|}/\Delta(a)$. Note that the asymptotic variance $\lim_{n \to \infty} \sigma_n^2 = \int_{R^d} c_\Xi^{(2)}(o, x) \, dx$ is finite and strictly positive iff $0 < \mathsf{E}|\Xi_0|^2 < \infty$. In this case we can find suitable upper and lower bounds of $c_\Xi^{(2)}(o, x) = \exp\{-\lambda \mathsf{E}|\Xi_0 \cup (\Xi_0 - x)|\} - \exp\{-2\lambda \mathsf{E}|\Xi_0|\}$ that lead to the inclusion

$$\exp\{-2\lambda \mathsf{E}|\Xi_0|\}(1 - \exp\{-\lambda \mathsf{E}|\Xi_0|\}) \frac{\mathsf{E}|\Xi_0|^2}{\mathsf{E}|\Xi_0|}$$

$$\leq \int_{R^d} c_\Xi^{(2)}(o, x) \, dx \leq \lambda \mathsf{E}|\Xi_0|^2 \exp\{-\lambda \mathsf{E}|\Xi_0|\}.$$

The estimate (4.2) enables us to apply to $\xi_n$ a well-known lemma on large deviations of a single random variable proved by Statulevičius (1966) which immediately provides the asymptotic relationships (1.11) and (1.12) as well as the Berry–Esseen bound (1.14) stated in Theorem 2. To be precise, according to the result by Statulevičius (1966), the relationships (1.11) and (1.12) are only valid in the narrower interval $0 \leq x \leq \delta^* \Delta_n$ for any $\delta^* < \delta_0(1+\delta_0)/2$, where $\delta_0 \in (0, 1)$ denotes the unique real root of $(1 - \delta)^3 = 6H_n \delta$. Indeed, since $H_n \geq 1/2$, by (3.14) for $n = 1$, we have $\delta_0(1 + \delta_0) < \delta_0/(1 - \delta_0)^3 = 1/6H_n \leq 1/(1 + 4H_n)$. Using again (4.1) and the inequality $\binom{k+l+1}{l} \leq 2^{k+l}$, we can estimate the coefficients (1.13) as

$$|\mu_k^{(n)}| \leq \frac{1}{(k+2)(k+3)} \sum_{l=1}^{k+1} 2^{k+l} \sum_{\substack{k_1+\cdots+k_l=k+1 \\ k_i \geq 1, i=1,\ldots,l}} \prod_{i=1}^{l} \left(\frac{\Delta(a)^{k_i} H(a)}{\sigma_n^2}\right)$$

$$= \frac{2^k \Delta(a)^{k+1}}{(k+2)(k+3)} \sum_{l=1}^{k+1} \binom{k}{l-1} \left(\frac{2H(a)}{\sigma_n^2}\right)^l$$

$$= \frac{4H_n \Delta(a)}{(k+2)(k+3)} (2\Delta(a)(1+4H_n))^k$$

for $k \geq 0$. Therefore, the series $\sum_{k \geq 0} \mu_k^{(n)} (x/\sigma_n \sqrt{|W_n|})^k$ converges absolutely for $|x| \leq \Delta_n/2(1+4H_n)$, and the $O$-terms in (1.11) and (1.12) can be easily verified by evaluating the remainder terms given by Statulevičius (1966). Thus, (1.11) and (1.12) are valid for the whole interval $0 \leq x \leq \Delta_n/2(1+4H_n)$, which completes the proof of Theorem 2. $\square$

The proof of the large deviations inequalities stated in the Theorem 3 relies on Chebychev's inequality combined with Lemma 7 and (3.14) [resp. Lemma 1 and (1.10)].

PROOF OF THEOREM 3. For any integer $N \geq 2$, the $N$th moment of a random variable $Y$ can be expressed by its cumulants $\Gamma_k(Y)$, $k = 1, \ldots, N$ [by inverting (2.4)] in the manner

$$(4.3) \qquad \mathsf{E} Y^N = \sum_{k=1}^{N} \frac{N!}{k!} \sum_{\substack{n_1+\cdots+n_k=N \\ n_i \geq 1, i=1,\ldots,k}} \frac{\Gamma_{n_1}(Y)}{n_1!} \cdots \frac{\Gamma_{n_k}(Y)}{n_k!}.$$

Consider (4.3) for $Y = |\Xi \cap W| - p|W|$ with $p = \mathsf{E}|\Xi \cap [0,1)^d|$. Since $\Gamma_1(Y) = \mathsf{E} Y = 0$ and, by Lemma 1 combined with (3.14), $|\Gamma_n(Y)| \leq c_{1,n-1}(\lambda)|W|$ for $n = 2, \ldots, N$, we are led to

$$|\mathsf{E} Y^N| \leq N! \sum_{k=1}^{[N/2]} \frac{|W|^k}{k!} \sum_{\substack{n_1+\cdots+n_k=N \\ n_i \geq 2, i=1,\ldots,k}} \frac{c_{1,n_1-1}(\lambda)}{n_1!} \cdots \frac{c_{1,n_k-1}(\lambda)}{n_k!}$$

$$\leq c_N^{(1)}(\lambda) \max\{|W|, |W|^{[N/2]}\},$$

where $c_N^{(1)}(\lambda)$ depends on the first $N$ moments of $|\Xi_0|$. Hence, for an even integer $s \geq 2$, (1.16) follows from Chebyshev's inequality. To prove (1.16) for any real $s \geq 2$ we next show

$$(4.4) \qquad \mathsf{E} ||\Xi \cap W| - p|W||^s \leq c_s^{(1)}(\lambda)|W|^{s/2} \qquad \text{for } s \geq 2$$

provided that $|W| \geq 1$. For this, we introduce the "truncated" stationary PGM $\Xi_W := \bigcup_{i \geq 1}(\Xi_i^W + X_i)$ generated by $\Pi_\lambda = \sum_{i \geq 1} \delta_{X_i}$ and the typical grain

$$\Xi_0^W := \Xi_0 \cap B_{R_0^W}(o),$$

where the random variable $R_0^W := \sup\{r > 0 : |\Xi_0 \cap B_r(o)| \leq |W|^\alpha\}$ takes the value $\infty$ if $|\Xi_0| < |W|^\alpha$. Here, we put $\alpha = N/2(N+2-s)$ if $N < s < N+2$ for some even integer $N \geq 2$. Define $Y_W := |\Xi_W \cap W| - |W|\mathsf{E}|\Xi_W \cap [0,1)^d|$ and let $\tilde{\Xi}_W$ denote the PGM with typical grain $\Xi_0 \setminus \Xi_0^W$. Since $\tilde{\Xi}_W \subset \Xi$ and $\Xi \setminus \Xi_W \subset \tilde{\Xi}_W$ we have

$$|Y| \leq |Y_W| + \max\{|\tilde{\Xi}_W \cap W|, |W|\mathsf{E}|\tilde{\Xi}_W \cap [0,1)^d|\},$$

which implies

(4.5) $\quad \mathsf{E}|Y|^s \leq 2^{s-1}\mathsf{E}|Y_W|^s + 2^{s-1}\mathsf{E}|\tilde{\Xi}_W \cap W|^s + 2^{s-1}(\lambda \mathsf{E}|\Xi_0 \setminus \Xi_0^W||W|)^s.$

By definition of $\Xi_0^W$,

$$\mathsf{E}|\Xi_0 \setminus \Xi_0^W|^k = \mathsf{E}|\Xi_0 \cap B_{R_0^W}^c(o)|^k \mathbf{1}_{\{|\Xi_0| \geq |W|^\alpha\}} \leq \mathsf{E}|\Xi_0|^s |W|^{-\alpha(s-k)}$$

$$\text{for } 0 < k \leq s.$$

Thus, $\mathsf{E}|\Xi_0 \setminus \Xi_0^W||W| \leq \mathsf{E}|\Xi_0|^s |W|^{1-\alpha(s-1)} \leq \mathsf{E}|\Xi_0|^s \sqrt{|W|}$. Next, applying (A.3),

$$\mathsf{E}|\tilde{\Xi}_W \cap W|^s \leq |W|^{s-N} \mathsf{E}\left(\sum_{i \geq 1} |((\Xi_i \setminus \Xi_i^W) + X_i) \cap W|\right)^N$$

$$\leq |W|^{s-N} N! \sum_{k=1}^{N} \frac{(\lambda|W|)^k}{k!} \sum_{\substack{n_1 + \cdots + n_k = N \\ n_i \geq 1, i=1,\ldots,k}} \prod_{i=1}^{k} \frac{\mathsf{E}|\Xi_0 \setminus \Xi_0^W|^{n_i}}{n_i!}$$

$$\leq c_1(N)|W|^{s-N} \max_{1 \leq k \leq N}\{(\lambda \mathsf{E}|\Xi_0|^s)^k |W|^{k(1-\alpha s) + \alpha N}\}$$

$$\leq c_2(N, \lambda)|W|^{s/2}.$$

Since, by Lyapunov's inequality, $\mathsf{E}|Y_W|^s \leq (\mathsf{E}|Y_W|^{N+2})^{s/(N+2)}$, we need only to verify that $\mathsf{E}|Y_W|^{N+2} \leq c_N(\lambda)|W|^{(N+2)/2}$, which in turn follows from (4.3) (with $Y_W$ and $N+2$ instead $Y$ and $N$) whenever $|\Gamma_{N+2}(|\Xi_W \cap W|)| \leq c_3(N,\lambda)|W|^{(N+2)/2}$. A thorough examination of the proofs of Lemmas 5–7 reveals that the constant $c_{1,n}(\lambda)$ in (3.14) takes on the form

$$c_{1,n}(\lambda) = \lambda \mathsf{E}|\Xi_0|^{n+1} + b_n^{(1)}(\lambda)\mathsf{E}|\Xi_0|^n + b_n^{(2)}(\lambda),$$

where $b_n^{(1)}(\lambda)$ and $b_n^{(2)}(\lambda)$ are given polynomials in $\exp\{\lambda \mathsf{E}|\Xi_0|\}$ and the first $n-1$ moments of $|\Xi_0|$. Hence, by $\mathsf{E}|\Xi_0^W|^{N+2} \leq |W|^{\alpha(N+2-s)}\mathsf{E}|\Xi_0|^s$, we get

the desired upper bound of $|\Gamma_{N+2}(|\Xi_W \cap W|)|$. Putting together the above estimates yields (4.4) and hence (1.16) is proved for any real $s \geq 2$.

To establish the second inequality (1.17), we use that $\Xi \cup (\Xi - x)$ is also a stationary PGM with typical grain $\Xi_0 \cup (\Xi_0 - x)$ and volume fraction $p(x) := \mathsf{P}(o \in \Xi \cup (\Xi - x))$. In view of the obvious decomposition

$$\hat{C}_W(x) - C(x) = \hat{p}_W - p + \frac{|(\Xi - x) \cap W|}{|W|} - p - \left(\frac{|(\Xi \cup (\Xi - x)) \cap W|}{|W|} - p(x)\right),$$

we obtain (1.17) by applying (1.16) to the three stationary PGMs $\Xi$, $\Xi - x$ and $\Xi \cup (\Xi - x)$. Finally, to prove the exponential inequality (1.18), we again employ a Chebyshev-type inequality. In this way we obtain, for $\varepsilon \geq 0$ and $0 \leq h \leq \rho/\Delta(a)$, that

$$\mathsf{P}(\hat{p}_W - p \geq \varepsilon) \leq \exp\{-h|W|(\varepsilon + p) + \log \mathsf{E} e^{h|\Xi \cap W|}\}$$
$$\leq \exp\left\{-h|W|\varepsilon + \frac{h^2 H(a)}{2}|W| \sum_{k \geq 2}(h\Delta(a))^{k-2}\right\}$$
$$\leq \exp\left\{-h|W|\varepsilon + \frac{h^2 H(a)}{2(1-\rho)}|W|\right\}.$$

Taking $h = \varepsilon(1-\rho)/H(a)$ for $0 \leq \varepsilon \leq H(a)\rho/\Delta(a)(1-\rho)$ proves the first part of (1.18), whereas the second part is obtained by setting $h = \rho/\Delta(a)$ in the latter inequality. $\square$

PROOF OF THEOREM 4. As in the proof of Theorem 3, using (1.2) and the notation $\hat{p}_n = \hat{p}_{W_n}$,

$$\mathsf{P}(\hat{p}_n - p \geq \varepsilon) \leq \exp\{|W_n|(L_n(h) - h(\varepsilon + p))\}$$

for any $h \geq 0$, whence, by virtue of Theorem 1, it follows that

$$\limsup_{n \to \infty} \frac{1}{|W_n|} \log \mathsf{P}(\hat{p}_n - p \geq \varepsilon) \leq g(h) \qquad \text{for } 0 \leq h < \frac{1}{\Delta(a)}.$$

Thus, the limit on the left-hand side is bounded from above by $\inf_{0 \leq h < 1/\Delta(a)} g(h)$.

Relationship (1.19) is proved as soon as we show that

$$(4.6) \qquad \liminf_{n \to \infty} \frac{1}{|W_n|} \log \mathsf{P}(\hat{p}_n - p \geq \varepsilon) \geq \inf_{0 \leq h < 1/\Delta(a)} g(h).$$

For brevity put $\zeta_n(\varepsilon) := |\Xi \cap W_n| - (\varepsilon + p)|W_n|$. Then, for any $\delta > 0$ and $h \geq 0$,

$$\mathsf{P}(\hat{p}_n - p \geq \varepsilon)$$
$$(4.7) \qquad \geq \mathsf{P}(\hat{p}_n - p \in (\varepsilon, \varepsilon + \delta])$$
$$\geq \exp\{|W_n|(L_n(h) - h(\varepsilon + p + \delta))\} \frac{\mathsf{E} \exp\{h\zeta_n(\varepsilon)\} \mathbf{1}_{\{\zeta_n(\varepsilon) \in (0, \delta|W_n|]\}}}{\mathsf{E} \exp\{h\zeta_n(\varepsilon)\}}.$$

Due to the properties of cumulant-generating functions [see, e.g., Dembo and Zeitouni (1998), page 27], the functions $L_n(\cdot)$ are convex on the whole real axis and $L_n''(h) > 0$ for every $h \in R^1$ (provided that $\sigma_n^2 > 0$). In view of Theorem 1, $L_n'(h) \underset{n\to\infty}{\longrightarrow} L'(h)$ and both $L_n'(h)$ (for $n \geq n_0$) and $L'(h)$ are strictly increasing for $0 \leq h < 1/\Delta(a)$. Likewise, $L_n''(h) \underset{n\to\infty}{\longrightarrow} L''(h)$ with $L''(h) \geq 0$ for $0 \leq h < 1/\Delta(a)$. Hence, for each $\varepsilon \in [0, \varepsilon^*)$ and sufficiently large $n$, there exists a unique $h_n = h_n(\varepsilon) \in [0, 1/\Delta(a))$ that satisfies the equation $L_n'(h_n) = \varepsilon + p$. Moreover, we have $h_n \underset{n\to\infty}{\longrightarrow} h_0$, where $h_0 = h_0(\varepsilon)$ is the unique solution of $L'(h) = \varepsilon + p$. Since $h \mapsto g(h)$ is a convex function and $g'(h_0) = 0$, it follows that $g(h_0) = \inf_{0 \leq h < 1/\Delta(a)} g(h)$. Consequently, putting $h = h_n$ on the right-hand side of (4.7) and taking into account that $L_n(h_n) \underset{n\to\infty}{\longrightarrow} L(h_0)$, we arrive at

$$
(4.8) \quad \begin{aligned}
\liminf_{n \to \infty} & \frac{\log \mathsf{P}(\hat{p}_n - p \geq \varepsilon)}{|W_n|} \\
& \geq g(h_0) - h_0 \delta + \liminf_{n \to \infty} \frac{1}{|W_n|} \log(G_n(\delta|W_n|) - G_n(0)),
\end{aligned}
$$

where the distribution function $G_n(x) = \mathsf{E} \exp\{h_n \zeta_n(\varepsilon)\} \mathbf{1}_{\{\zeta_n(\varepsilon) \leq x\}} / \mathsf{E} \exp\{h_n \zeta_n(\varepsilon)\}$ possesses the Fourier–Stieltjes transform $\hat{G}_n(t) = \mathsf{E} \exp\{(it + h_n)\zeta_n(\varepsilon)\} / \mathsf{E} \exp\{h_n \zeta_n(\varepsilon)\}$. Using (1.2) and $L_n'(h_n) = \varepsilon + p$ we can write

$$
\begin{aligned}
\log \hat{G}_n(t) &= |W_n|(-it(\varepsilon + p) + L_n(it + h_n) - L_n(h_n)) \\
&= -|W_n| t^2 \int_0^1 (1 - \vartheta) L_n''(i\vartheta t + h_n) \, d\vartheta \qquad \text{if } |t| + h_n < 1/\Delta(a),
\end{aligned}
$$

where the last line is obtained by partial integration of $L_n''(it\vartheta + h_n)$ with respect to $\vartheta \in [0, 1]$. An application of Theorem 1 shows that $\log \hat{G}_n(t/\sqrt{|W_n|}) \underset{n\to\infty}{\longrightarrow} -t^2 \times L''(h_0)/2$ for all $t \in R^1$, which in turn implies $G_n(x\sqrt{|W_n|}) \underset{n\to\infty}{\longrightarrow} \Phi(x/\sqrt{L''(h_0)})$ provided that $L''(h_0) > 0$. In this case, $G_n(\delta|W_n|) - G_n(0) \underset{n\to\infty}{\longrightarrow} 1/2$, proving (4.6) and, thus, the desired relationship (1.19) holds. If $L''(h_0(\varepsilon_0)) = 0$ for certain $\varepsilon_0 \in (0, \varepsilon^*)$, then there exists some $\eta > 0$ such that $L''(h_0(\varepsilon)) > 0$ for $\varepsilon \in [\varepsilon_0 - \eta, \varepsilon_0 + \eta] \setminus \{\varepsilon_0\}$. Since $\log \mathsf{P}(\hat{p}_n - p \geq \varepsilon)$ is nonincreasing in $\varepsilon$, it follows that

$$
\begin{aligned}
g(h_0(\varepsilon_0 + \eta)) &\leq \liminf_{n \to \infty} \frac{\log \mathsf{P}(\hat{p}_n - p \geq \varepsilon_0)}{|W_n|} \\
&\leq \limsup_{n \to \infty} \frac{\log \mathsf{P}(\hat{p}_n - p \geq \varepsilon_0)}{|W_n|} \\
&\leq g(h_0(\varepsilon_0 - \eta)).
\end{aligned}
$$

Hence, having in mind the continuity of $h_0(\cdot)$, the proof of Theorem 4 is finished.

□

## APPENDIX

For any random closed set $\Xi$ (defined as $(\mathfrak{A}, \sigma_f)$-measurable mapping $\Xi : [\Omega, \mathfrak{A}, \mathsf{P}] \mapsto [\mathcal{F}, \sigma_f]$, where $\mathcal{F}$ is the (metrizable) space of all closed set in $R^d$ and $\sigma_f$ is its Matheron $\sigma$-field), the mapping $(x, \omega) \mapsto \mathbf{1}_{\Xi(\omega)}(x)$ is $(\mathfrak{B}(R^d) \otimes \mathfrak{A})$-measurable since $(x, F) \mapsto \mathbf{1}_F(x)$ is $(\mathfrak{B}(R^d) \otimes \sigma_f)$-measurable [see Matheron (1975), Chapter 2]. As we see below, the PGM (1.1) is no longer P-a.s. closed in $R^d$ if $\mathsf{E}|\Xi_0| < \infty$ and $\mathsf{E}|\Xi_0 + B_\varepsilon(o)| = \infty$ for any $\varepsilon > 0$.

To preserve the $(\mathfrak{B}(R^d) \otimes \mathfrak{A})$-measurability of the indicator function $\mathbf{1}_{\Xi(\omega)}(x)$ (needed to apply Fubini's theorem), we define the not necessarily closed PGM (1.1) as a set-valued, measurable mapping of an independently marked Poisson process $\Pi_{\lambda,Q} := \sum_{i \geq 1} \delta_{[X_i, \Xi_i]}$ with mark distribution $Q(\mathcal{L}) = \mathsf{P}(\Xi_0 \in \mathcal{L})$, $\mathcal{L} \in \mathfrak{B}(\mathcal{K})$.

More precisely, let $M_\mathcal{K}$ denote the space of all integer-valued measures $\psi$ on $[R^d \times \mathcal{K}, \mathfrak{B}(R^d) \otimes \mathfrak{B}(\mathcal{K})]$ that satisfy $\psi(B \times \mathcal{K}) < \infty$ for each bounded $B \in \mathfrak{B}(R^d)$ and let $\mathfrak{M}_\mathcal{K}$ be the $\sigma$-field generated by the sets $\{\psi \in M_\mathcal{K} : \psi(B \times \mathcal{L}) = n\}$ for $n \geq 0$, bounded $B \in \mathfrak{B}(R^d)$ and $\mathcal{L} \in \mathfrak{B}(\mathcal{K})$. Each $\psi \in M_\mathcal{K}$ admits a representation $\psi = \sum_{i \geq 1} \delta_{[x_i(\psi), k_i(\psi)]}$ as a sum of Dirac measures with respect to the at most countable set of atoms $[x_i(\psi), k_i(\psi)]$, $i \geq 1$, where each atom is counted according to its multiplicity. Note that the mappings $M_\mathcal{K} \ni \psi \mapsto [x_i(\psi), k_i(\psi)] \in R^d \times \mathcal{K}$ are measurable [see Matthes, Kerstan and Mecke (1978)]. Finally, define

(A.1) $$\Xi(\psi) := \bigcup_{i \geq 1} (x_i(\psi) + k_i(\psi)).$$

PROPOSITION 1. *The mapping $(x, \psi) \mapsto \xi(x, \psi) := \mathbf{1}_{\Xi(\psi)}(x)$ is $(\mathfrak{B}(R^d) \otimes \mathfrak{M}_\mathcal{K})$-measurable, that is, $\{(x, \psi) \in R^d \times M_\mathcal{K} : x \notin \Xi(\psi)\} \in \mathfrak{B}(R^d) \otimes \mathfrak{M}_\mathcal{K}$.*

PROOF. Let $\mathcal{B}$ be the countable set of open balls in $R^d$ that have rational radii and midpoints with rational coordinates. For any sequence $\{K_n, n \geq 1\}$ in $\mathcal{K}$ that satisfies $K_n \uparrow R^d$ put $\Xi_n(\psi) := \bigcup_{i : x_i \in K_n}(x_i(\psi) + k_i(\psi))$. Obviously, $\Xi_n(\psi) \in \mathcal{K}$ and $\Xi(\psi) = \bigcup_{n \geq 1} \Xi_n(\psi)$.

It can be readily checked that the set $\xi^{-1}(\{0\}) = \{(x, \psi) \in R^d \times M_\mathcal{K} : x \notin \Xi(\psi)\}$ coincides with $\bigcap_{n \geq 1} \bigcup_{B \in \mathcal{B}} (B \times \{\psi \in M_\mathcal{K} : \Xi_n(\psi) \cap B = \varnothing\})$. However, this set belongs to $\mathfrak{B}(R^d) \times \mathfrak{M}_\mathcal{K}$ if $\{\psi : \Xi_n(\psi) \cap B = \varnothing\} = \{\psi : \sum_{i \geq 1} \mathbf{1}_{\mathcal{R}_B \cap (K_n \times \mathcal{K})}(x_i(\psi), k_i(\psi)) = 0\} \in \mathfrak{M}_\mathcal{K}$ for any $n \geq 1$ and $B \in \mathcal{B}$, where $\mathcal{R}_B := \{(x, K) \in R^d \times \mathcal{K} : (x + K) \cap B \neq \varnothing\}$. Since the mapping $\psi \mapsto \sum_{i \geq 1} f(x_i(\psi), k_i(\psi))$ is $\mathfrak{M}_\mathcal{K}$-measurable whenever $f: R^d \times \mathcal{K} \mapsto R^1$ is $(\mathfrak{B}(R^d) \otimes \mathfrak{B}(\mathcal{K}))$-measurable [see Matthes, Kerstan and Mecke (1978)], we need only to verify that $\mathcal{R}_B \in \mathfrak{B}(R^d) \otimes \mathfrak{B}(\mathcal{K})$ for $B \in \mathcal{B}$. Since there exists a sequence of closed balls $B_n$ such that $B_n \uparrow B$ and thus $\mathcal{R}_{B_n} \uparrow \mathcal{R}_B$, it suffices to show $\mathcal{R}_C \in \mathfrak{B}(R^d) \otimes \mathfrak{B}(\mathcal{K})$ for any closed ball $C$. For this, remember that $\mathcal{K}^D := \{K \in \mathcal{K} : K \cap D = \varnothing\} \in \mathfrak{B}(\mathcal{K})$ for any open ball $D$ [see Matheron (1975)]. Thus, the proof of Proposition 1 is completed by noting that

$$\mathcal{R}_C^c = \{(x, K) \in R^d \times \mathcal{K} : K \cap (C - x) = \varnothing\} = \bigcup_{B \in \mathcal{B}} (B \times \mathcal{K}^{D(B,C)}),$$

where $D(B, C) = \{u - v : u \in C, v \in B\}$ is an open ball for $B \in \mathcal{B}$. $\square$

Let there be given an unmarked point process $\Psi = \sum_{i \geq 1} \delta_{X_i}$ on $R^d$ and a random compact set $\Xi_0$ with distribution $Q = \mathsf{P} \circ \Xi_0^{-1}$ on $[\mathcal{K}, \mathfrak{B}(\mathcal{K})]$. The corresponding independently marked point process $\Psi_Q = \sum_{i \geq 1} \delta_{[X_i, \Xi_i]}$ on $R^d$ with mark distribution $Q$ is then defined [see Daley and Vere-Jones (1988) or Stoyan, Kendall and Mecke (1995)] to be a random element (over $[\Omega, \mathfrak{A}, \mathsf{P}]$) that takes values in $[M_\mathcal{K}, \mathfrak{M}_\mathcal{K}]$, the distribution of which is uniquely determined by its probability generating functional $G_{\Psi_Q}[v] = \mathsf{E}(\prod_{i \geq 1} v(X_i, \Xi_i)) = G_\Psi[v_Q]$ (see Section 2), where $v_Q(\cdot) := \int_\mathcal{K} v(\cdot, K) Q(dK)$ and the function $v : R^d \times \mathcal{K} \mapsto [0, 1]$ is Borel measurable such that $1 - v(\cdot, K)$ has bounded support for all $K \in \mathcal{K}$. In the special case of a stationary independently marked Poisson process $\Pi_{\lambda,Q} = \sum_{i \geq 1} \delta_{[X_i, \Xi_i]}$, the shape of $G_{\Pi_\lambda}[\cdot]$ implies that

$$(A.2) \qquad G_{\Pi_{\lambda,Q}}[v] = \exp\left\{\lambda \int_{R^d} \int_\mathcal{K} (v(x, K) - 1) Q(dK) \, dx\right\}$$

and, furthermore, the Campbell-type formula

$$(A.3) \quad \mathsf{E}\left(\sum_{i_1,\ldots,i_k \geq 1}^* \prod_{j=1}^k f_j(X_{i_j}, \Xi_{i_j})\right) = \lambda^k \prod_{j=1}^k \int_{R^d} \int_\mathcal{K} f_j(x, K) Q(dK) \, dx$$

holds for any measurable functions $f_1, \ldots, f_k : R^d \times \mathcal{K} \mapsto [0, \infty]$, where the sum $\sum^*$ on the left-hand side of (A.3) stretches over $k$-tuples of pairwise distinct indices.

As announced, we conclude the Appendix by showing that, under the assumption $\mathsf{E}|\Xi_0| < \infty$, the condition $\mathsf{E}|\Xi_0 + B_\varepsilon(o)| < \infty$ for some $\varepsilon > 0$ is not only sufficient as shown by Heinrich (1992), but even necessary for the closedness of the stationary PGM $\Xi = \Xi(\Pi_{\lambda,Q})$.

PROPOSITION 2. *Let $\Xi_0$ be a compact typical grain of the PGM (1.1) that satisfies $\mathsf{E}|\Xi_0| < \infty$ and $\mathsf{E}|\Xi_0 + B_\varepsilon(o)| = \infty$ for any $\varepsilon > 0$. Then $\mathsf{P}(\Xi$ is closed in $R^d) = 0$.*

PROOF. Choose $K_n \in \mathcal{K}$, $n \geq 1$, such that $K_n \uparrow R^d$ as $n \to \infty$ and let $\Xi_n(\psi)$ be defined as in the proof of Proposition 1. Obviously, $\{\psi \in M_\mathcal{K} : \Xi_n(\psi) \cap B_\varepsilon(o) = \varnothing\} \downarrow \{\psi \in M_\mathcal{K} : \Xi(\psi) \cap B_\varepsilon(o) = \varnothing\}$ as $n \to \infty$. Furthermore, since $(\Xi_i + X_i) \cap B_\varepsilon(o) \neq \varnothing$ iff $-X_i \in \Xi_i + B_\varepsilon(o)$, we find, using (A.2), that

$$\mathsf{P}(\Xi_n(\Pi_{\lambda,Q}) \cap B_\varepsilon(o) = \varnothing) = \mathsf{E} \prod_{i \geq 1}(1 - \mathbf{1}_{K_n}(X_i)\mathbf{1}_{\Xi_i + B_\varepsilon(o)}(-X_i))$$

$$= \exp\{-\lambda \mathsf{E}|(-K_n) \cap (\Xi_0 + B_\varepsilon(o))|\}.$$

By the monotone convergence theorem and our assumptions,

$$\lim_{n \to \infty} \mathsf{E}|(-K_n) \cap (\Xi_0 + B_\varepsilon(o))| = \mathsf{E}|\Xi_0 + B_\varepsilon(o)| = \infty,$$

which means for the stationary PGM $\Xi = \Xi(\Pi_{\lambda,Q})$ that

$$\mathsf{P}(\Xi \cap B_\varepsilon(o) \neq \varnothing) = 1 - \lim_{n \to \infty} \mathsf{P}(\Xi_n(\Pi_{\lambda,Q}) \cap B_\varepsilon(o) = \varnothing) = 1 \qquad \text{for any } \varepsilon > 0.$$

Thus,

$$\mathsf{P}(\Xi \text{ is closed}) = \mathsf{P}\bigg(\Xi \text{ is closed}, \bigcap_{m \geq 1}\{\Xi \cap B_{1/m}(o) \neq \varnothing\}\bigg)$$

$$\leq \mathsf{P}(o \in \Xi) = 1 - \exp\{-\lambda \mathsf{E}|\Xi_0|\} < 1.$$

Similarly, $\mathsf{P}(\Xi \text{ is closed}) \leq \mathsf{P}(x_1 \in \Xi, \ldots, x_n \in \Xi)$ for any $x_1, \ldots, x_n \in R^d$. In view of (1.4)–(1.6), the probability $p_\Xi^{(n)}(x_1, \ldots, x_n)$ is arbitrarily close to $(\mathsf{P}(o \in \Xi))^n$ whenever the distances between the points $x_1, \ldots, x_n$ are sufficiently large. This proves the assertion of Proposition 2. $\square$

**Acknowledgment.** I thank one of the referees for helpful comments and suggestions.

## REFERENCES

DALEY, D. J. and VERE-JONES, D. (1988). *An Introduction to the Theory of Point Processes.* Springer, New York. MR950166

DEMBO, A. and ZEITOUNI, O. (1998). *Large Deviations Techniques and Applications*, 2nd ed. Springer, New York. MR1619036

FRITZ, J. (1970). Generalization of McMillan's theorem to random set functions. *Studia Sci. Math. Hungar.* **5** 369–394. MR293956

GEORGII, H.-O. and ZESSIN, H. (1993). Large deviations and maximum entropy principle for marked point random fields. *Probab. Theory Related Fields* **96** 177–204. MR1227031

GÖTZE, F., HEINRICH, L. and HIPP, C. (1995). $m$-dependent random fields with analytic cumulant generating function. *Scand. J. Statist.* **22** 183–195. MR1339750

GREENBERG, W. (1971). Thermodynamic states of classical systems. *Comm. Math. Phys.* **22** 259–268. MR351345

HALL, P. (1988). *Introduction to the Theory of Coverage Processes.* Wiley, New York. MR973404

Heinrich, L. (1992). On existence and mixing properties of germ-grain models. *Statistics* **23** 271–286. MR1237805

Heinrich, L. and Molchanov, I. S. (1999). Central limit theorem for a class of random measures associated with germ-grain models. *Adv. in Appl. Probab.* **31** 283–314. MR1724553

Heinrich, L. and Schmidt, V. (1985). Normal convergence of multidimensional shot noise and rates of this convergence. *Adv. in Appl. Probab.* **17** 709–730. MR809427

Ivanoff, G. (1982). Central limit theorems for point processes. *Stochastic Process. Appl.* **12** 171–186. MR651902

Mase, S. (1982). Asymptotic properties of stereological estimators of volume fraction for stationary random sets. *J. Appl. Probab.* **19** 111–126. MR644424

Matheron, G. (1975). *Random Sets and Integral Geometry*. Wiley, New York. MR385969

Matthes, K., Kerstan, J. and Mecke, J. (1978). *Infinitely Divisible Point Processes*. Wiley, Chichester. MR517931

Molchanov, I. S. (1997). *Statistics of the Boolean Model for Practitioners and Mathematicians*. Wiley, Chichester.

Piau, D. (1999). Large deviations and Young measures for a Poissonian model of biphased material. *Ann. Appl. Probab.* **9** 706–718. MR1722279

Rota, G.-C. (1964). On the foundations of combinatorial theory: I. Theory of Möbius functions. *Z. Wahrsch. Verw. Gebiete* **2** 340–368. MR174487

Ruelle, D. (1964). Cluster properties of correlation functions of classical gases. *Rev. Mod. Phys.* **36** 580–586. MR167259

Ruelle, D. (1969). *Statistical Mechanics*: *Rigorous Results*. Benjamin, Amsterdam. MR289084

Saulis, L. and Statulevičius, V. (1991). *Limit Theorems for Large Deviations*. Kluwer Academic, Dordrecht. MR1171883

Statulevičius, V. A. (1966). On large deviations. *Z. Wahrsch. Verw. Gebiete* **6** 133–144. MR221560

Stoyan, D., Kendall, W. S. and Mecke, J. (1995). *Stochastic Geometry and Its Applications*, 2nd ed. Wiley, Chichester. MR895588

Institute of Mathematics
University of Augsburg
D-86135 Augsburg
Germany
e-mail: heinrich@math.uni-augsburg.de